\documentclass[preprint,12pt,1p]{elsarticle}
\makeatletter
\def\ps@pprintTitle{%
  \let\@oddhead\@empty
  \let\@evenhead\@empty
  \let\@oddfoot\@empty
  \let\@evenfoot\@oddfoot
}
\makeatother
\usepackage{graphics}
\usepackage{epsfig}
\usepackage{mathptmx}
\usepackage{amsmath}
\usepackage{amssymb}
\usepackage{hyperref}
\usepackage{fullpage}
\usepackage{amsthm}
\usepackage{geometry}
\theoremstyle{plain}
\newtheorem{theorem}{Theorem}[section]

\theoremstyle{definition}

\theoremstyle{example}

\theoremstyle{remark}

\newlength{\defbaselineskip}
\journal{Computer Aided Geometric Design}
\numberwithin{equation}{section}
\begin{document}
\begin{frontmatter}
\title{Some hypergeometric summation theorems and reduction formulas via Laplace transform method}
\author{M. I. Qureshi}
\author{*Showkat Ahmad Dar}
\ead{showkat34@gmail.com}
\address{miqureshi\_delhi@yahoo.co.in,~showkat34@gmail.com}
\address{Department of Applied Sciences and Humanities , \\Faculty of  Engineering and Technology, \\Jamia Millia Islamia ( Central University), New Delhi, 110025, India.}
\cortext[cor1]{Corresponding author}
\begin{abstract}
In this paper, we obtain analytical solutions of Laplace transform based some generalized class of the hyperbolic integrals in terms of hypergeometric functions ${}_{3}F_{2}(\pm1)$, ${}_{4}F_{3}(\pm1)$, ${}_{5}F_{4}(\pm1)$,\\${}_{6}F_{5}(\pm1)$, ${}_{7}F_{6}(\pm1)$ and ${}_{8}F_{7}(\pm1)$ with suitable convergence conditions, by using some algebraic properties of Pochhammer symbols. In addition, reduction formulas for ${}_{4}F_{3}(1)$, ${}_{7}F_{6}(-1)$ and some new summation theorems (not recorded earlier in the literature of hypergeometric functions) for ${}_{3}F_{2}(-1), {}_{6}F_{5}(\pm1)$, ${}_{7}F_{6}(\pm1)$ and ${}_{8}F_{7}(\pm1)$ are obtained.
 \\
 \\
\textit{2010 AMS Classification: 33C05; 33C20; 44A10; 33B15  }
 \end{abstract}
\end{frontmatter}
\section{Introduction and Preliminaries}
For the sake of conciseness of this paper, we use the following notations \\
$~~~~~~~~~~~\mathbb{N}:=\{1,2,...\};~~~~~~\mathbb{N}_{0}:=\mathbb{N}\cup\{0\};~~~~~~\mathbb{Z}_{0}^{-}:=\mathbb{Z}^{-}\cup\{0\}=\{0,-1,-2,-3,...\},$\\
where the symbols $\mathbb{N}$ and $\mathbb{Z}$ denote the set of natural numbers and integers; as usual, the symbols $\mathbb{R}$ and $\mathbb{C}$ denote the set of real and complex numbers respectively. \\
Here the notation $(\lambda)_{\upsilon}~(\lambda, \upsilon\in\mathbb{C})$ denotes the Pochhammer's symbol (or the shifted factorial,
 since $(1)_{n}=n! )$ is defined, in general, by
\begin{equation}\label{GRI04}
(\lambda)_{\upsilon}:=\frac{\Gamma(\lambda+\upsilon)}{\Gamma(\lambda)}=\begin{cases} 1 , \quad~~~~~~~~~~~~~ (\upsilon=0~;~\lambda\in\mathbb{C}\backslash\{0\}) \\ \lambda(\lambda+1)...(\lambda+n-1),  \quad (\upsilon=n\in\mathbb{N}~;~\lambda\in\mathbb{C}). \\
 \end{cases}
 \end{equation}
A natural generalization of Gauss hypergeometric series ${}_{2}F_{1}$ is the general hypergeometric series ${}_{p}F_{q}$  
 with $p$ numerator parameters $\alpha_{1},  ... , \alpha_{p}$ and $q$ denominator parameters $\beta_{1}, ..., \beta_{q}$. It is defined by\\
\begin{equation}\label{GRI05}
{}_{p}F_{q}\left(\ \begin{array}{lll}\alpha_{1},...,\alpha_{p}~;~\\\beta_{1}, ..., \beta_{q}~;~\end{array} z\right)
=\sum_{n=0}^{\infty}\frac{(\alpha_{1})_{n}  ...  (\alpha_{p})_{n}}{(\beta_{1})_{n} ... (\beta_{q})_{n}}\frac{z^{n}}{n!}~,\newline
\end{equation}
where $\alpha_{i}\in\mathbb{C}~(i=1,...,p)$ and $\beta_{j}\in\mathbb{C}\setminus \mathbb{Z}_{0}^{-}~(j=1,...,q)~\left(\ \mathbb{Z}_{0}^{-}:=\{0,-1,-2,...\}\right)$ and \\
$\left(\ p,~q\in\mathbb{N}_{0}:=\mathbb{N}\cup\{0\}=\{0,1,2,...\}\right)$.
The ${}_{p}F_{q}(\cdot)$ series in eq.(\ref{GRI05}) is convergent for $|z|<\infty$ if $p\leq q$, and for $|z|<1$ if $p=q+1$.
Furthermore, if we set
\begin{equation}\label{GRI06}
\omega=\left(\ \sum_{j=1}^{q}\beta_{j}-\sum_{i=1}^{p}\alpha_{i}\right),\newline
\end{equation}
it is known that the ${}_{p}F_{q}$ series, with $p=q+1$, is\\
(i) absolutely convergent for $|z|=1$ if $\Re(\omega)>0,$\\
(ii) conditionally convergent for $|z|=1, z\neq1$, if $-1<Re(\omega)\leq0$.\\
The  binomial function is given by
\begin{equation}\label{GRI08}
(1-z)^{-a}={}_{1}F_{0}\left(\begin{array}{lll}a~;\\  \overline{~~~};\end{array} z \right)
=\sum_{n=0}^{\infty}\frac{(a)_{n}}{n!}z^{n},
\end{equation}
~~~~~~~~~where $|z|<1,~~a\in\mathbb{C}$.\\
Next we collect some results that we will need in the sequel.\\
\begin{equation}
B_{z}(\alpha,\beta)=\int_{0}^{z}t^{\alpha-1}(1-t)^{\beta-1}dt,~~~~~0<z\leq1,
\end{equation}
\begin{equation}
{}_{2}F_{1}\left(\begin{array}{lll}a,~~b;\\  1+b,;\end{array} z\right)=\frac{b}{z^{b}}B_{z}(b,1-a),
\end{equation}
where $B_{z}$ is incomplete beta function.\\
The Dixon's theorem for ${}_{3}F_{2}$ with positive unit argument is given by \\
\begin{equation}\label{GRI11}
{}_{3}F_{2}\left(\begin{array}{lll}a,~~b,~~c~~~~~~~~~~~~~~~~~~;\\  1+a-b,~1+a-c;\end{array} 1\right)
=\frac{\Gamma(1+a-b)\Gamma(1+a-c)\Gamma(1+\frac{a}{2})\Gamma(1+\frac{a}{2}-b-c)}{\Gamma(1+\frac{a}{2}-b)\Gamma(1+\frac{a}{2}-c)\Gamma(1+a)\Gamma(1+a-b-c)},
\end{equation}
~~~~~~~~where $\Re(a-2b-2c)>-2;~ 1+a-b,~1+a-c\in\mathbb{C}\backslash \mathbb{Z}_{0}^{-}$,\\
and when $c=1+\frac{a}{2}$ in the eq.(\ref{GRI11}), we get
\begin{equation}
{}_{3}F_{2}\left(\begin{array}{lll}a,~~1+\frac{a}{2},~~ b;\\ \frac{a}{2},~~1+a-b;~\end{array} 1\right)=0,
\end{equation}
~~~~~~~~where $\Re(b)<0;~\frac{a}{2}, 1+a-b\in\mathbb{C}\backslash \mathbb{Z}_{0}^{-}$. \\
The following summation theorems ${}_{3}F_{2}(-1)$ is given by
\begin{equation}\label{FGG}
{}_{3}F_{2}\left(\begin{array}{lll}a,~~1+\frac{a}{2},~~ b;\\ \frac{a}{2},~~1+a-b;~\end{array} -1\right)
=\frac{\Gamma(1+a-b)\Gamma(\frac{1+a}{2})}{\Gamma(\frac{1+a}{2}-b)\Gamma(1+a)},
\end{equation}
~~~~~~~~where $\Re(b)<\frac{1}{2}$; ~$\frac{a}{2}, 1+a-b\in\mathbb{C}\backslash \mathbb{Z}_{0}^{-}$.\\
The eq.(\ref{FGG}) can be obtained by setting $2d=1+a$ and $c=b$ in the eq.(\ref{GRI12}).\\
Contiguous function relations \cite[p.71, Q.N0.21, part(13)]{R}
\begin{equation}\label{GGG}
(\beta-\gamma+1)~{}_{2}F_{1}\left(\begin{array}{lll}\alpha,~\beta;\\ \gamma~~~~;~\end{array} z\right)=\beta~{}_{2}F_{1}\left(\begin{array}{lll}\alpha,~\beta+1;\\ \gamma~~~~~~~~~~~~;~\end{array} z\right)
-(\gamma-1)~{}_{2}F_{1}\left(\begin{array}{lll}\alpha,~\beta;\\ \gamma-1~;~\end{array} z\right).
\end{equation}
In the both sides of eq.(\ref{GGG}) replace $z$ by $zt$, multiply by $t^{h-1}(1-t)^{g-h-1}$, integrate with respect to $t$ over the interval $(0,1)$ and using the definition of beta function, after simplification we get
\begin{equation}\label{GGH}
(\beta-\gamma+1)~{}_{3}F_{2}\left(\begin{array}{lll}\alpha,~\beta,~h;\\ \gamma,~g~~~~~;~\end{array} z\right)=\beta~{}_{3}F_{2}\left(\begin{array}{lll}\alpha,~\beta+1,~h;\\ \gamma,~g~~~~~~~~~~~;~\end{array} z\right)-(\gamma-1)~{}_{3}F_{2}\left(\begin{array}{lll}\alpha,~\beta,~h;\\ \gamma-1,~g;~\end{array} z\right).
\end{equation}
In the eq.(\ref{GGH}) put $\alpha=a,~\beta=c,~h=b,~\gamma=1+d,~g=1+c$, after simplification we get
\begin{equation}\label{GGR}
{}_{3}F_{2}\left(\begin{array}{lll}a,~b,~c~~~~~~~;\\ 1+c,~1+d;~\end{array} z\right)=\frac{c}{(c-d)}~{}_{2}F_{1}\left(\begin{array}{lll}a,~b;\\ d+1;~\end{array} z\right)-\frac{d}{(c-d)}~{}_{3}F_{2}\left(\begin{array}{lll}a,~b,~c;\\ c+1,d;~\end{array} z\right).
\end{equation}
Put $d=b$ in eq.(\ref{GGR}), we get\\
When $c\neq b$ and $1+b,1+c\in\mathbb{C}\backslash \mathbb{Z}_{0}^{-}$, then
\begin{equation}\label{FRI}
{}_{3}F_{2}\left(\begin{array}{lll}a,~b,~c~~~~~~~;\\ 1+b,~1+c;~\end{array} z\right)=\left(\frac{c}{c-b}\right){}_{2}F_{1}\left(\begin{array}{lll}a,~b;\\ 1+b;\end{array} z\right)-\left(\frac{b}{c-b}\right){}_{2}F_{1}\left(\begin{array}{lll}a,~c;\\ 1+c;\end{array} z\right).
\end{equation}
When $z=1$ in the eq.(\ref{FRI}) and using Gauss classical summation theorem, we get
\begin{equation}
  {}_{3}F_{2}\left(\begin{array}{lll}a,~b,~c~~~~~~~;\\ 1+b,~1+c;~\end{array} 1\right)=\frac{bc}{(c-b)}\Gamma(1-a)\left\{\frac{\Gamma(b)}{\Gamma(1+b-a)}-\frac{\Gamma(c)}{\Gamma(1+c-a)}\right\},
\end{equation}
where $c\neq b; \Re(a)<1$; $1+b,1+c\in\mathbb{C}\backslash \mathbb{Z}_{0}^{-}$,\\
\\
The classical summation theorems for hypergeometric series ${}_{4}F_{3}(\pm1)$ \cite[ p.28, eq.(4.4.3)]{B2} are given by
\begin{equation}\label{GRI12}
{}_{4}F_{3}\left(\begin{array}{lll}a,~~1+\frac{a}{2},~~ c,~~d~~~~~~~~~~~~~;\\ \frac{a}{2},~~1+a-c,~~ 1+a-d;~\end{array} -1\right)
=\frac{\Gamma(1+a-c)\Gamma(1+a-d)}{\Gamma(1+a)\Gamma(1+a-c-d)}\newline,
\end{equation}
~~~~~~~~provided $\Re(a-2c-2d)>-2 ;~ \frac{a}{2}, 1+a-c, 1+a-d\in\mathbb{C}\backslash \mathbb{Z}_{0}^{-}$.\\
When $e=(1+a)/2$ in the eq.(\ref{GRI0}), we get
\begin{multline}\label{GRI102}
{}_{4}F_{3}\left(\begin{array}{lll}a,~~1+\frac{a}{2},~~ c,~~d~~~~~~~~~~~~~;\\ \frac{a}{2},~~1+a-c,~~ 1+a-d;~\end{array} 1\right)\\
=\frac{\Gamma(1+a-c)\Gamma(1+a-d)\Gamma(\frac{1+a}{2})\Gamma(\frac{1+a}{2}-c-d)}{\Gamma(1+a)\Gamma(\frac{1+a}{2}-d)\Gamma(\frac{1+a}{2}-c)\Gamma(1+a-c-d)},\newline
\end{multline}
~~~~~~~~provided $\Re(2c+2d-a)<1;~ \frac{a}{2}, 1+a-c, 1+a-d\in\mathbb{C}\backslash \mathbb{Z}_{0}^{-}$.\\
Another classical summation theorem for hypergeometric series ${}_{5}F_{4}(1)$ \cite[p.27,eq.(4.4.1)]{B2}, is given by
\begin{multline}\label{GRI0}
{}_{5}F_{4}\left(\begin{array}{lll}a,~~1+\frac{a}{2},~~ c,~~d,~~e~~~~~~~~~~~~~~~~~~~~~~~~~;\\ \frac{a}{2},~~1+a-c,~~ 1+a-d,~~ 1+a-e;~\end{array} 1\right)\\
=\frac{\Gamma(1+a-c)\Gamma(1+a-d)\Gamma(1+a-e)\Gamma(1+a-c-d-e)}{\Gamma(1+a)\Gamma(1+a-d-e)\Gamma(1+a-c-e)\Gamma(1+a-c-d)},
\end{multline}
provided $\Re(a-c-d-e)>-1 ;~ \frac{a}{2}, 1+a-c, 1+a-d,1+a-e\in\mathbb{C}\backslash \mathbb{Z}_{0}^{-}$.\\
Laplace transform of any constant $k$ is given by
\begin{equation}\label{GRI22}
\mathcal{L}[k;q]=\int_{0}^{\infty}e^{-qt}k~dt=\frac{k}{q},
\end{equation}
~~~~~provided
\begin{equation}\label{GRI23}
\Re(q)>0.
\end{equation}
The Digamma function (or Psi function)  is given by
\begin{equation*}
  \Psi(x)=\frac{d}{dx}\{\ln\Gamma(x)\}=\frac{\Gamma~^{'}(x)}{\Gamma(x)},
\end{equation*}
\begin{equation}
  =-\gamma+(x-1)\sum_{n=0}^{\infty}\frac{1}{(n+1)(n+x)},
\end{equation}
where $\gamma(=0.57721566490...)$ being the Euler-Mascheroni constant
and $\Psi^{'}(x)=\frac{d}{dx}\{\Psi(x)\}$ is called trigamma function,\\
where
\begin{equation}
\Psi^{'}(x)=\sum_{k=0}^{\infty}\frac{1}{(x+k)^{2}}=\frac{1}{x^{2}}~{}_{3}F_{2}\left(\begin{array}{lll} 1,~x,x~~~;\\ 1+x,~1+x;\end{array} 1\right),
\end{equation}
OR
\begin{equation}
{}_{3}F_{2}\left(\begin{array}{lll} 1,~x,x~~~;\\ 1+x,~1+x;\end{array} 1\right)=x^{2}\Psi^{'}(x),
\end{equation}
and
\begin{equation}\label{GRI24}
{}_{3}F_{2}\left(\begin{array}{lll} 1,~a,b~~~;\\ 1+a,~1+b;\end{array} 1\right)=\frac{ab}{(b-a)}\left[\Psi(b)-\Psi(a)\right],~~~~b\neq a.
\end{equation}
Properties of Digamma function \cite{G,G3}
\begin{equation}\label{GRI25}
\Psi(1+x)=\Psi(x)+\frac{1}{x},
\end{equation}
\begin{equation}\label{GRI262}
\Psi(1-x)=\Psi(x)+\pi\cot(\pi x),
\end{equation}
where $x\neq0,\pm1,\pm2,\pm3,...$.
\begin{equation}
\Psi\left(\frac{1}{2}+x\right)-\Psi\left(\frac{1}{2}-x\right)=\pi \tan(\pi x),
\end{equation}
 where $x\neq\pm\frac{1}{2},\pm\frac{3}{2},\pm\frac{5}{2},...$.\\
Lower case beta function of one variable which is related with Digamma function, is given by
\begin{equation}\label{GRI25}
\beta(x)=\sum_{k=0}^{\infty}\frac{(-1)^{k}}{(k+x)}=\frac{1}{x}~{}_{2}F_{1}\left(\begin{array}{lll} 1,~x;\\ 1+x;\end{array} -1\right),
\end{equation}
\begin{equation}\label{GRI26}
  ~~~~~~~~~=\frac{1}{2}\bigg[\Psi\left(\frac{1+x}{2}\right)-\Psi\left(\frac{x}{2}\right)\bigg],~~~~~~~~~x\in\mathbb{C}\backslash \mathbb{Z}_{0}^{-},
\end{equation}
and
\begin{equation}\label{GRI000}
{}_{3}F_{2}\left(\begin{array}{lll} 1,~a,b~~~;\\ 1+a,~1+b;\end{array} -1\right)=\frac{ab}{(b-a)}\left[\beta(a)-\beta(b)\right],~~~~b\neq a,
\end{equation}
\begin{equation}
{}_{3}F_{2}\left(\begin{array}{lll} 1,~x,x~~~;\\ 1+x,~1+x;\end{array} -1\right)=-x^{2}\frac{d}{dx}\{\beta(x)\}=-x^{2}\beta^{'}(x).
\end{equation}
Hypergeometric forms  of some trigonometric ratios \cite[pp.137-138, eq.(4.2.6); eq.(4.2.9); eq.(4.2.3)]{hy} and \cite[pp.48-49, eqns(20,23,26)]{QI} , are given by
\begin{equation}\label{GRI31}
 \tan(z)=\frac{8z}{(\pi^{2}-4z^{2})}~~{}_{3}F_{2}\left(\begin{array}{lll} 1,~\frac{1}{2}+\frac{z}{\pi},~~\frac{1}{2}-\frac{z}{\pi};\\  \frac{3}{2}+\frac{z}{\pi},~~\frac{3}{2}-\frac{z}{\pi}~~~~;\end{array} 1\right),
 \end{equation}
 where $z\in\mathbb{C}\setminus\{\pm\frac{\pi}{2},~\pm\frac{3\pi}{2},~\pm\frac{5\pi}{2},...\}$,
 \begin{equation}\label{GRI32}
 \sec(z)=\frac{4\pi}{(\pi^{2}-4z^{2})}~~{}_{4}F_{3}\left(\begin{array}{lll} 1,~\frac{3}{2},~~\frac{1}{2}+\frac{z}{\pi},~~\frac{1}{2}-\frac{z}{\pi};\\  \frac{1}{2},~~\frac{3}{2}+\frac{z}{\pi},~~\frac{3}{2}-\frac{z}{\pi}~~~~;\end{array} -1\right),
 \end{equation}
 where $z\in\mathbb{C}\setminus\{\pm\frac{\pi}{2},~\pm\frac{3\pi}{2},~\pm\frac{5\pi}{2},...\}$,
 \begin{eqnarray}\label{GRI33}
 \sec^{2}(z)=\frac{4}{(2z-\pi)^{2}}~~{}_{3}F_{2}\left(\begin{array}{lll} 1,~\frac{1}{2}-\frac{z}{\pi},~~\frac{1}{2}-\frac{z}{\pi};\\  \frac{3}{2}-\frac{z}{\pi},~~\frac{3}{2}-\frac{z}{\pi}~;\end{array} 1\right)+\\
 +\frac{4}{(2z+\pi)^{2}}~~{}_{3}F_{2}\left(\begin{array}{lll} 1,~\frac{1}{2}+\frac{z}{\pi},~\frac{1}{2}+\frac{z}{\pi};\\  \frac{3}{2}+\frac{z}{\pi},~\frac{3}{2}+\frac{z}{\pi}~;\end{array} 1\right),
 \end{eqnarray}
 where $z\in\mathbb{C}\setminus\{\pm\frac{\pi}{2},~\pm\frac{3\pi}{2},~\pm\frac{5\pi}{2},...\}$.\\
 \\
$~~~~~~$ The plan of this paper is as follows. First, we obtain  generalized class and analytical solutions of some hyperbolic integrals in terms of  ${}_{6}F_{5}(\pm1)$, ${}_{7}F_{6}(\pm1)$, ${}_{8}F_{7}(\pm1)$ shown in \texttt{section 4}. Special class of some hyperbolic integrals in terms of ${}_{3}F_{2}(\pm1),{}_{4}F_{3}(\pm1),{}_{5}F_{4}(\pm1)$ are given in \texttt{section 5}.  We apply suitable product formulas associated with hyperbolic function in special class of hyperbolic integrals  given in \texttt{section 6}.
 Moreover, we find two reduction formulas and some new summation theorems given in \texttt{sections 2-3} by comparing the similar integrals.
\section{Some new summation theorems}
\begin{theorem}
The first summation theorem ${}_{6}F_{5}(-1)$ holds true:
 \begin{multline}
 {}_{6}F_{5}\left(\begin{array}{lll}v,~~1+\frac{v}{2},~~\frac{v}{2}-\frac{a}{2c}-\frac{b}{2c},~~\frac{v}{2}+\frac{a}{2c}+\frac{b}{2c},
 ~~\frac{v}{2}-\frac{a}{2c}+\frac{b}{2c},~~\frac{v}{2}+\frac{a}{2c}-\frac{b}{2c}~~~~~~~~;\\  \frac{v}{2},~1+\frac{v}{2}-\frac{a}{2c}-\frac{b}{2c},~1+\frac{v}{2}+\frac{a}{2c}+\frac{b}{2c},~1+\frac{v}{2}-\frac{a}{2c}+\frac{b}{2c},~1+\frac{v}{2}+\frac{a}{2c}-\frac{b}{2c},~~;\end{array} -1\right)\\
 =\frac{\{(vc)^{2}-(a+b)^{2}\}\{(vc)^{2}-(a-b)^{2}\}}{16vabc^{2}\Gamma(v)}\bigg[\Gamma\left(\frac{vc+a+b}{2c}\right)\Gamma\left(\frac{vc-a-b}{2c}\right)\\
   -\Gamma\left(\frac{vc+a-b}{2c}\right)\Gamma\left(\frac{vc-a+b}{2c}\right)\bigg],
 \end{multline}
 \end{theorem}
 where $\Re(v)<4,~\Re(c)>0,~\Re(vc\pm a\pm b)>0;~\frac{v}{2},~1+\frac{v}{2}\pm\frac{a}{2c}\pm\frac{b}{2c}\in\mathbb{C}\backslash \mathbb{Z}_{0}^{-}$.\\
  \textbf{Proof}: Comparing the two equations (\ref{GRI45}) and (\ref{sss1}),we get a summation theorem  for ${}_{6}F_{5}(-1)$.
 \begin{theorem}
 The second summation theorem ${}_{6}F_{5}(1)$ holds true:
 \begin{multline}
 {}_{6}F_{5}\left(\begin{array}{lll}v,~~1+\frac{v}{2},~~\frac{v}{2}-\frac{a}{2c}-\frac{b}{2c},~~\frac{v}{2}+\frac{a}{2c}+\frac{b}{2c},
 ~~\frac{v}{2}-\frac{a}{2c}+\frac{b}{2c},~~\frac{v}{2}+\frac{a}{2c}-\frac{b}{2c}~~~~~~~~;\\  \frac{v}{2},~1+\frac{v}{2}-\frac{a}{2c}-\frac{b}{2c},~1+\frac{v}{2}+\frac{a}{2c}+\frac{b}{2c},~1+\frac{v}{2}-\frac{a}{2c}+\frac{b}{2c},~1+\frac{v}{2}+\frac{a}{2c}-\frac{b}{2c},~~;\end{array} 1\right)\\
 =\frac{\{(vc)^{2}-(a+b)^{2}\}\{(vc)^{2}-(a-b)^{2}\}}{16vabc^{2}\Gamma(v)}
 \bigg[\Gamma\left(\frac{vc+a+b}{2c}\right)\Gamma\left(\frac{vc-a-b}{2c}\right)\frac{\cos\left(\frac{(a+b)\pi}{2c}\right)}{\cos(\frac{v\pi}{2})}\\
   -\Gamma\left(\frac{vc+a-b}{2c}\right)\Gamma\left(\frac{vc-a+b}{2c}\right)\frac{\cos\left(\frac{(a-b)\pi}{2c}\right)}{\cos(\frac{v\pi}{2})}\bigg],
 \end{multline}
 \end{theorem}
  where $\Re(v)<3,~\Re(c)>0,~\Re(vc\pm a\pm b)>0;~\frac{v}{2},~1+\frac{v}{2}\pm\frac{a}{2c}\pm\frac{b}{2c}\in\mathbb{C}\backslash \mathbb{Z}_{0}^{-}$.\\
  \textbf{Proof}: Comparing the two equations (\ref{GRI46}) and (\ref{ss2}),we get a summation theorem  for ${}_{6}F_{5}(1)$.
=

 \begin{theorem}
  The  third summation theorem ${}_{6}F_{5}(-1)$ holds true:
 \begin{multline}
 {}_{6}F_{5}\left(\begin{array}{lll}1,\frac{3}{2},\frac{1}{2}-\frac{a}{2b}-\frac{c}{2b},
 \frac{1}{2}-\frac{a}{2b}+\frac{c}{2b},\frac{1}{2}+\frac{a}{2b}+\frac{c}{2b},\frac{1}{2}+\frac{a}{2b}-\frac{c}{2b}~;\\  \frac{1}{2},\frac{3}{2}-\frac{a}{2b}-\frac{c}{2b}
 ,\frac{3}{2}-\frac{a}{2b}+\frac{c}{2b},\frac{3}{2}+\frac{a}{2b}+\frac{c}{2b},\frac{3}{2}+\frac{a}{2b}-\frac{c}{2b};\end{array} -1\right)\\
 =\frac{\pi(b-a-c)(b+a+c)(b-a+c)(b+a-c)}{2^{2}acb^{2}}
 \frac{\sin\left(\frac{a\pi}{2b}\right)\sin\left(\frac{c\pi}{2b}\right)}{\{\cos\left(\frac{c\pi}{b}\right)+\cos\left(\frac{a\pi}{b}\right)\}},
 \end{multline}
  \end{theorem}
  where $\Re(b)>0,~\Re(b\pm a\pm c)>0,~\frac{3}{2}\pm\frac{a}{2b}\pm\frac{c}{2b}\in\mathbb{C}\backslash \mathbb{Z}_{0}^{-}$,$\frac{a\pm c}{b}\in\mathbb{C}\backslash\{\pm1,\pm3,\pm5,...\}$.\\
  \textbf{Proof}: Comparing the two equations (\ref{FAF81}) and (\ref{FAF811}), we get a summation theorem  for ${}_{6}F_{5}(-1)$.
  \begin{theorem}
   The fourth summation theorem ${}_{7}F_{6}(1)$ holds true:
 \begin{multline}\label{FA56}
 {}_{7}F_{6}\left(\begin{array}{lll}v,1+\frac{v}{2}-\frac{\sqrt{a^{2}-b^{2}}}{2c},1+\frac{v}{2}+\frac{\sqrt{a^{2}-b^{2}}}{2c},\frac{v}{2}-\frac{a}{2c}-\frac{b}{2c},
 \frac{v}{2}-\frac{a}{2c}+\frac{b}{2c},\frac{v}{2}+\frac{a}{2c}+\frac{b}{2c},\frac{v}{2}+\frac{a}{2c}-\frac{b}{2c}~;\\  \frac{v}{2}-\frac{\sqrt{a^{2}-b^{2}}}{2c},\frac{v}{2}+\frac{\sqrt{a^{2}-b^{2}}}{2c},1+\frac{v}{2}-\frac{a}{2c}-\frac{b}{2c}
 ,1+\frac{v}{2}-\frac{a}{2c}+\frac{b}{2c},1+\frac{v}{2}+\frac{a}{2c}+\frac{b}{2c},1+\frac{v}{2}+\frac{a}{2c}-\frac{b}{2c};\end{array} 1\right)\\
 =\frac{\{(vc)^{2}-(a+b)^{2}\}\{(vc)^{2}-(a-b)^{2}\}}{16(v^{2}ac^{3}-a^{3}c+ab^{2}c)\Gamma(v)}
 \bigg[\Gamma\left(\frac{vc+a+b}{2c}\right)\Gamma\left(\frac{vc-a-b}{2c}\right)\frac{\sin\left(\frac{(a+b)\pi}{2c}\right)}{\sin(\frac{v\pi}{2})}\\
   +\Gamma\left(\frac{vc+a-b}{2c}\right)\Gamma\left(\frac{vc-a+b}{2c}\right)\frac{\sin\left(\frac{(a-b)\pi}{2c}\right)}{\sin(\frac{v\pi}{2})}\bigg],
 \end{multline}
 \end{theorem}
 where  $\Re(v)<2,~\Re(c)>0,~\Re(vc\pm a\pm b)>0;~\frac{v}{2}\pm\frac{\sqrt{a^{2}-b^{2}}}{2c},~1+\frac{v}{2}\pm\frac{a}{2c}\pm\frac{b}{2c}\in\mathbb{C}\backslash \mathbb{Z}_{0}^{-}$.\\
  \textbf{Proof}: Comparing the two equations (\ref{GRI48}) and (\ref{ss3}),we get a summation theorem  for ${}_{7}F_{6}(1)$.
  \begin{theorem}
  The fifth summation theorem ${}_{7}F_{6}(1)$ holds true:
 \begin{multline}\label{FA56}
 {}_{7}F_{6}\left(\begin{array}{lll}1,\frac{3}{2}-\frac{\sqrt{a^{2}-c^{2}}}{2b},\frac{3}{2}+\frac{\sqrt{a^{2}-c^{2}}}{2b},\frac{1}{2}-\frac{a}{2b}-\frac{c}{2b},
 \frac{1}{2}-\frac{a}{2b}+\frac{c}{2b},\frac{1}{2}+\frac{a}{2b}+\frac{c}{2b},\frac{1}{2}+\frac{a}{2b}-\frac{c}{2b}~;\\  \frac{1}{2}-\frac{\sqrt{a^{2}-c^{2}}}{2b},\frac{1}{2}+\frac{\sqrt{a^{2}-c^{2}}}{2b},\frac{3}{2}-\frac{a}{2b}-\frac{c}{2b}
 ,\frac{3}{2}-\frac{a}{2b}+\frac{c}{2b},\frac{3}{2}+\frac{a}{2b}+\frac{c}{2b},\frac{3}{2}+\frac{a}{2b}-\frac{c}{2b};\end{array} 1\right)\\
 =\frac{\pi(b-a-c)(b+a+c)(b-a+c)(b+a-c)}{4(ab^{3}-a^{3}b+abc^{2})}
 \frac{\sin\left(\frac{a\pi}{b}\right)}{\{\cos\left(\frac{c\pi}{b}\right)+\cos\left(\frac{a\pi}{b}\right)\}},
 \end{multline}
 \end{theorem}
  where $\Re(b)>0,~\Re(b\pm a\pm c)>0;~\frac{1}{2}\pm\frac{\sqrt{a^{2}-c^{2}}}{2b},~\frac{3}{2}\pm\frac{a}{2b}\pm\frac{c}{2b}\in\mathbb{C}\backslash \mathbb{Z}_{0}^{-}$,~$\frac{a\pm c}{b}\in\mathbb{C}\backslash\{\pm1,\pm3,\pm5,...\}$.\\
 \textbf{Proof}: Comparing the two equations (\ref{FAF71}) and (\ref{FAF710}), we get a summation theorem  for ${}_{7}F_{6}(1)$.
  \begin{theorem}
 The sixth summation theorem ${}_{7}F_{6}(-1)$ holds true:
  \begin{multline}
 {}_{7}F_{6}\left(\begin{array}{lll}1,\frac{3}{2}-\frac{\sqrt{a^{2}-c^{2}}}{2b},\frac{3}{2}+\frac{\sqrt{a^{2}-c^{2}}}{2b},\frac{1}{2}-\frac{a}{2b}-\frac{c}{2b},
 \frac{1}{2}-\frac{a}{2b}+\frac{c}{2b},\frac{1}{2}+\frac{a}{2b}+\frac{c}{2b},\frac{1}{2}+\frac{a}{2b}-\frac{c}{2b}~;\\  \frac{1}{2}-\frac{\sqrt{a^{2}-c^{2}}}{2b},\frac{1}{2}+\frac{\sqrt{a^{2}-c^{2}}}{2b},\frac{3}{2}-\frac{a}{2b}-\frac{c}{2b}
 ,\frac{3}{2}-\frac{a}{2b}+\frac{c}{2b},\frac{3}{2}+\frac{a}{2b}+\frac{c}{2b},\frac{3}{2}+\frac{a}{2b}-\frac{c}{2b};\end{array} -1\right)\\
  =\frac{\pi(b-a-c)(b+a+c)(b-a+c)(b+a-c)}{2(ab^{3}-a^{3}b+abc^{2})}
  \frac{\cos\left(\frac{a\pi}{2b}\right)\cos\left(\frac{c\pi}{2b}\right)}{\{\cos\left(\frac{c\pi}{b}\right)+\cos\left(\frac{a\pi}{b}\right)\}}\\
     -\frac{(b-a-c)(b+a+c)(b-a+c)(b+a-c)}{4(ab^{3}-a^{3}b+abc^{2})}\left[\beta\left(\frac{(a+b+c)}{2b}\right)+\beta\left(\frac{(a+b-c)}{2b}\right)\right],
 \end{multline}
  \end{theorem}
   where $\Re(b)>0,~\Re(b\pm a\pm c)>0;~\frac{1}{2}\pm\frac{\sqrt{a^{2}-c^{2}}}{2b},~\frac{3}{2}\pm\frac{a}{2b}\pm\frac{c}{2b}\in\mathbb{C}\backslash \mathbb{Z}_{0}^{-}$,$\frac{a\pm c}{b}\in\mathbb{C}\backslash\{\pm1,\pm3,\pm5,...\}$.\\
  \textbf{Proof}:Comparing the two equations (\ref{FAF101}) and (\ref{FAF111}),we get a summation theorem for ${}_{7}F_{6}(-1)$.
   \begin{theorem}
   The seventh summation theorem ${}_{8}F_{7}(-1)$ holds true:
  \begin{multline}
 {}_{8}F_{7}\left(\begin{array}{lll}v,1+\frac{v}{2},1+\frac{v}{2}-\frac{\sqrt{a^{2}+b^{2}}}{2c},1+\frac{v}{2}+\frac{\sqrt{a^{2}+b^{2}}}{2c},\frac{v}{2}-\frac{a}{2c}-\frac{b}{2c},
 \frac{v}{2}-\frac{a}{2c}+\frac{b}{2c},\frac{v}{2}+\frac{a}{2c}+\frac{b}{2c},\frac{v}{2}+\frac{a}{2c}-\frac{b}{2c}~;\\
  \frac{v}{2}, \frac{v}{2}-\frac{\sqrt{a^{2}+b^{2}}}{2c},\frac{v}{2}+\frac{\sqrt{a^{2}+b^{2}}}{2c},1+\frac{v}{2}-\frac{a}{2c}-\frac{b}{2c}
 ,1+\frac{v}{2}-\frac{a}{2c}+\frac{b}{2c},1+\frac{v}{2}+\frac{a}{2c}+\frac{b}{2c},1+\frac{v}{2}+\frac{a}{2c}-\frac{b}{2c};\end{array} -1\right)\\
 =\frac{\{(vc)^{2}-(a+b)^{2}\}\{(vc)^{2}-(a-b)^{2}\}}{8(v^{3}c^{4}-a^{2}vc^{2}-b^{2}vc^{2})\Gamma(v)}
 \bigg[\Gamma\left(\frac{vc+a+b}{2c}\right)\Gamma\left(\frac{vc-a-b}{2c}\right)\\
   +\Gamma\left(\frac{vc+a-b}{2c}\right)\Gamma\left(\frac{vc-a+b}{2c}\right)\bigg],
 \end{multline}
 \end{theorem}
 where $\Re(v)<2,~\Re(c)>0,~\Re(vc\pm a\pm b)>0;\frac{v}{2},~\frac{v}{2}\pm\frac{\sqrt{a^{2}+b^{2}}}{2c},~1+\frac{v}{2}\pm\frac{a}{2c}\pm\frac{b}{2c}\in\mathbb{C}\backslash \mathbb{Z}_{0}^{-}$.\\
  \textbf{Proof}: Comparing the two equations (\ref{GRI49}) and (\ref{ss4}),we get a summation theorem  for ${}_{8}F_{7}(-1)$.
 \begin{theorem}
   The eigth summation theorem ${}_{8}F_{7}(1)$ holds true:
  \begin{multline}
 {}_{8}F_{7}\left(\begin{array}{lll}v,1+\frac{v}{2},1+\frac{v}{2}-\frac{\sqrt{a^{2}+b^{2}}}{2c},1+\frac{v}{2}+\frac{\sqrt{a^{2}+b^{2}}}{2c},\frac{v}{2}-\frac{a}{2c}-\frac{b}{2c},
 \frac{v}{2}-\frac{a}{2c}+\frac{b}{2c},\frac{v}{2}+\frac{a}{2c}+\frac{b}{2c},\frac{v}{2}+\frac{a}{2c}-\frac{b}{2c}~;\\
  \frac{v}{2}, \frac{v}{2}-\frac{\sqrt{a^{2}+b^{2}}}{2c},\frac{v}{2}+\frac{\sqrt{a^{2}+b^{2}}}{2c},1+\frac{v}{2}-\frac{a}{2c}-\frac{b}{2c}
 ,1+\frac{v}{2}-\frac{a}{2c}+\frac{b}{2c},1+\frac{v}{2}+\frac{a}{2c}+\frac{b}{2c},1+\frac{v}{2}+\frac{a}{2c}-\frac{b}{2c};\end{array} 1\right)\\
 =\frac{\{(vc)^{2}-(a+b)^{2}\}\{(vc)^{2}-(a-b)^{2}\}}{8(v^{3}c^{4}-a^{2}vc^{2}-b^{2}vc^{2})\Gamma(v)}
 \bigg[\Gamma\left(\frac{vc+a+b}{2c}\right)\Gamma\left(\frac{vc-a-b}{2c}\right)\frac{\cos\left(\frac{(a+b)\pi}{2c}\right)}{\cos(\frac{v\pi}{2})}\\
   +\Gamma\left(\frac{vc+a-b}{2c}\right)\Gamma\left(\frac{vc-a+b}{2c}\right)\frac{\cos\left(\frac{(a-b)\pi}{2c}\right)}{\cos(\frac{v\pi}{2})}\bigg],
 \end{multline}
 \end{theorem}
 where $\\Re(v)<1,~Re(c)>0,~\Re(vc\pm a\pm b)>0;~\frac{v}{2},~\frac{v}{2}\pm\frac{\sqrt{a^{2}+b^{2}}}{2c},~1+\frac{v}{2}\pm\frac{a}{2c}\pm\frac{b}{2c}\in\mathbb{C}\backslash \mathbb{Z}_{0}^{-}$.\\
  \textbf{Proof}: Comparing the two equations (\ref{GRI50}) and (\ref{ss5}),we get a summation theorem  for ${}_{8}F_{7}(1)$.
   \begin{theorem}
   The ninth summation theorem ${}_{8}F_{7}(-1)$ holds true:
  \begin{multline}
 {}_{8}F_{7}\left(\begin{array}{lll}1,\frac{3}{2},\frac{3}{2}-\frac{\sqrt{a^{2}+c^{2}}}{2b},\frac{3}{2}+\frac{\sqrt{a^{2}+c^{2}}}{2b},\frac{1}{2}-\frac{a}{2b}-\frac{c}{2b},
 \frac{1}{2}-\frac{a}{2b}+\frac{c}{2b},\frac{1}{2}+\frac{a}{2b}+\frac{c}{2b},\frac{1}{2}+\frac{a}{2b}-\frac{c}{2b}~;\\  \frac{1}{2}, \frac{1}{2}-\frac{\sqrt{a^{2}+c^{2}}}{2b},\frac{1}{2}+\frac{\sqrt{a^{2}+c^{2}}}{2b},\frac{3}{2}-\frac{a}{2b}-\frac{c}{2b}
 ,\frac{3}{2}-\frac{a}{2b}+\frac{c}{2b},\frac{3}{2}+\frac{a}{2b}+\frac{c}{2b},\frac{3}{2}+\frac{a}{2b}-\frac{c}{2b};\end{array} -1\right)\\
 =\frac{\pi(b-a-c)(b+a+c)(b-a+c)(b+a-c)}{2(b^{4}-a^{2}b^{2}-c^{2}b^{2})}
 \frac{\cos\left(\frac{a\pi}{2b}\right)\cos\left(\frac{c\pi}{2b}\right)}{\{\cos\left(\frac{c\pi}{b}\right)+\cos\left(\frac{a\pi}{b}\right)\}},
 \end{multline}
 \end{theorem}
 where $\Re(b)>0,~\Re(b\pm a\pm c)>0;~\frac{1}{2}\pm\frac{\sqrt{a^{2}+c^{2}}}{2b},~\frac{3}{2}\pm\frac{a}{2b}\pm\frac{c}{2b}\in\mathbb{C}\backslash \mathbb{Z}_{0}^{-}$,$\frac{a\pm c}{b}\in\mathbb{C}\backslash\{\pm1,\pm3,\pm5,...\}$.\\
  \textbf{Proof}: Comparing the two equations (\ref{FAF09}) and (\ref{FAF109}),we get a summation theorem for ${}_{8}F_{7}(-1)$.
\begin{theorem}
The tenth  summation theorem ${}_{3}F_{2}(-1)$ holds true:
 \begin{equation}\label{FA55}
{}_{3}F_{2}\left(\begin{array}{lll}1,~\frac{1}{2}-\frac{a}{2b},~~\frac{1}{2}+\frac{a}{2b};\\  \frac{3}{2}-\frac{a}{2b},~~\frac{3}{2}+\frac{a}{2b};\end{array} -1\right)=\frac{(b^{2}-a^{2})}{2ab}\bigg[\frac{\pi}{2}\sec\left(\frac{\pi a}{2b}\right)-\beta\left(\frac{a+b}{2b}\right)\bigg],
 \end{equation}
 \begin{equation}\label{FA77}
=\frac{(b^{2}-a^{2})}{8ab}\left[\Psi\left(\frac{3b-a}{4b}\right)-\Psi\left(\frac{b-a}{4b}\right)-\Psi\left(\frac{3b+a}{4b}\right)+\Psi\left(\frac{b+a}{4b}\right)\right],
 \end{equation}
 where $\Re(b)>0,~\Re(b\pm a)>0,~\frac{3}{2}\pm\frac{a}{2b}\in\mathbb{C}\backslash \mathbb{Z}_{0}^{-}$,~ $\frac{a}{b}\in\mathbb{C}\backslash\{\pm1,\pm3,\pm5,...\}$.
 \end{theorem}
  \textbf{Proof}: The summation formula with unit negative argument (\ref{FA55}) is obtained by comparing the two solutions  of the integral $\int_{0}^{\infty}\frac{\sinh(ax)}{\cosh(b x)}dx$,  given in (\ref{FA53}), (\ref{FA54}) and is not available in the literature of the hypergeometric summation theorem. The above eq.(\ref{FA77}) is obtained by using properties of beta function of one variable (\ref{GRI26}) and (\ref{GRI000}).
  \section{ Some reduction formulas}
\begin{theorem}
The first reduction formula holds true:
\begin{equation} \label{FAF17}
{}_{4}F_{3}\left(\begin{array}{lll}2,2,~\frac{1}{2}+\frac{a}{2b},~\frac{1}{2}-\frac{a}{2b};\\  1,~\frac{5}{2}+\frac{a}{2b},~\frac{5}{2}-\frac{a}{2b}~~~;\end{array} 1\right)
=\frac{(9b^{2}-a^{2})}{8b^{2}}~~{}_{3}F_{2}\left(\begin{array}{lll}1,~\frac{1}{2}-\frac{a}{2b},~~\frac{1}{2}+\frac{a}{2b};\\  \frac{3}{2}-\frac{a}{2b},~~\frac{3}{2}+\frac{a}{2b}~~~~;\end{array} 1\right),
\end{equation}
\end{theorem}
where $\Re(b)>0,~\Re(b\pm a)>0,~\frac{3}{2}\pm\frac{a}{2b}\in\mathbb{C}\backslash\mathbb{Z}_{0}^{-}$;~$\frac{a}{b}\in\mathbb{C}\backslash\{\pm1,\pm3,\pm5,...\}$.\\
\textbf{Proof}: The reduction formula (\ref{FAF17}) is obtained by comparing the two integrals (\ref{FA1}) and (\ref{FA33}).
 \begin{theorem}
The second reduction formula holds true:
 \begin{eqnarray*}
 {}_{7}F_{6}\left(\begin{array}{lll}v,1+\frac{v}{2}-\frac{\sqrt{a^{2}-b^{2}}}{2c},1+\frac{v}{2}+\frac{\sqrt{a^{2}-b^{2}}}{2c},
 \frac{v}{2}-\frac{a}{2c}-\frac{b}{2c},\frac{v}{2}+\frac{a}{2c}+\frac{b}{2c},\frac{v}{2}-\frac{a}{2c}+\frac{b}{2c},\frac{v}{2}+\frac{a}{2c}-\frac{b}{2c};\\
 \frac{v}{2}-\frac{\sqrt{a^{2}-b^{2}}}{2c},\frac{v}{2}+\frac{\sqrt{a^{2}-b^{2}}}{2c},1+\frac{v}{2}-\frac{a}{2c}-\frac{b}{2c},
 1+\frac{v}{2}+\frac{a}{2c}+\frac{b}{2c},1+\frac{v}{2}-\frac{a}{2c}+\frac{b}{2c},1+\frac{v}{2}+\frac{a}{2c}-\frac{b}{2c};\end{array} -1\right)
 \end{eqnarray*}
 \begin{multline}\label{FAF18}
 =\frac{\{(vc)^{2}-(a+b)^{2}\}\{(vc)^{2}-(a-b)^{2}\}}{4(v^{2}ac^{2}-a^{3}+ab^{2})}\times\\\times
 \bigg[\frac{1}{(vc-a-b)}~{}_{2}F_{1}\left(\begin{array}{lll}v,\frac{vc-a-b}{2c}~~~;\\ 1+\frac{vc-a-b}{2c};\end{array} -1\right)
 -\frac{1}{(vc+a+b)}~{}_{2}F_{1}\left(\begin{array}{lll}v,\frac{vc+a+b}{2c}~~~;\\ 1+\frac{vc+a+b}{2c};\end{array} -1\right)\\
 +\frac{1}{(vc-a+b)}~{}_{2}F_{1}\left(\begin{array}{lll}v,\frac{vc-a+b}{2c}~~~;\\ 1+\frac{vc-a+b}{2c};\end{array} -1\right)
 -\frac{1}{(vc+a-b)}~{}_{2}F_{1}\left(\begin{array}{lll}v,\frac{vc+a-b}{2c}~~~;\\ 1+\frac{vc+a-b}{2c};\end{array} -1\right)\bigg],
 \end{multline}
 where $\Re(c)>0,~\Re(vc\pm a\pm b)>0;~\frac{v}{2}\pm\frac{\sqrt{a^{2}-b^{2}}}{2c},~1+\frac{v}{2}\pm\frac{a}{2c}\pm\frac{b}{2c}\in\mathbb{C}\backslash \mathbb{Z}_{0}^{-}$.
 \end{theorem}
 \textbf{Proof}: Comparing the integral (\ref{GRI47}) and its companion (\ref{cc9}), we get the reduction formula (\ref{FAF18}).
 The integral representation of ${}_{2}F_{1}(-1)$ type hypergeometric functions involved in reduction formula (\ref{FAF18}) is given below
\begin{equation}\label{FRS}
{}_{2}F_{1}\left(\begin{array}{lll}a,~b,;\\ 1+b;\end{array} -1\right)=b\int_{0}^{1}\frac{t^{b-1}}{(1+t)^{a}}dt,~~~~~~\Re(b)>0.
\end{equation}
The definite integral (\ref{FRS}) can be solved by using suitable numerical methods (for example composite trapezoidal rule, composite Simpson's $1/3$ rule, composite Simpson's $3/8$ rule, composite Boole rule, composite Midris two rules, composite Weddle rule, composite Sadiq rule,
Gauss-Legendre three points formula, Gauss-Chebyshev three points formula, Radau three points formula, Lobatto three points formula etc ).
  \section{Generalized class of some hyperbolic integrals in terms of  ${}_{6}F_{5}(\pm1)$, ${}_{7}F_{6}(\pm1)$ and ${}_{8}F_{7}(\pm1)$}
 Many authors have studied some definite integrals containing the integrands as a quotient
of hyperbolic functions.  Mainly,  V. H. Moll et.al evaluated some definite integrals given in the table of Gradshteyn and Ryzhik \cite{G,G3}, by using the change of independent variables. We obtain generalizations and analytical solutions of some hyperbolic integrals, using  hypergeometric approach and Laplace transform method.
\begin{theorem}
The first generalized hyperbolic integral holds true:
 \begin{multline}\label{GRI45}
 \int_{0}^{\infty}\frac{\sinh(ax)\sinh(bx)}{\cosh^{v}(cx)}dx=\frac{2^{v+1}v abc}{(vc-a-b)(vc+a+b)(vc-a+b)(vc+a-b)}\times\\\times
 ~~{}_{6}F_{5}\left(\begin{array}{lll}v,~~1+\frac{v}{2},~~\frac{v}{2}-\frac{a}{2c}-\frac{b}{2c},~~\frac{v}{2}+\frac{a}{2c}+\frac{b}{2c},
 ~~\frac{v}{2}-\frac{a}{2c}+\frac{b}{2c},~~\frac{v}{2}+\frac{a}{2c}-\frac{b}{2c}~~~~~~~~~~~~~~~;\\  \frac{v}{2},~1+\frac{v}{2}-\frac{a}{2c}-\frac{b}{2c},~1+\frac{v}{2}+\frac{a}{2c}+\frac{b}{2c},~1+\frac{v}{2}-\frac{a}{2c}+\frac{b}{2c},~1+\frac{v}{2}+\frac{a}{2c}-\frac{b}{2c},~~;\end{array} -1\right),
 \end{multline}
 \end{theorem}
where $\Re(v)<4,~\Re(c)>0,~\Re(vc\pm a\pm b)>0;~\frac{v}{2},~1+\frac{v}{2}\pm\frac{a}{2c}\pm\frac{b}{2c}\in\mathbb{C}\backslash \mathbb{Z}_{0}^{-}$.
 \begin{theorem}
 The second generalized hyperbolic integral holds true:
 \begin{multline}\label{GRI46}
 \int_{0}^{\infty}\frac{\sinh(ax)\sinh(bx)}{\sinh^{v}(cx)}dx=\frac{2^{v+1}v abc}{(vc-a-b)(vc+a+b)(vc-a+b)(vc+a-b)}\times\\\times
 ~~{}_{6}F_{5}\left(\begin{array}{lll}v,~~1+\frac{v}{2},~~\frac{v}{2}-\frac{a}{2c}-\frac{b}{2c},~~\frac{v}{2}+\frac{a}{2c}+\frac{b}{2c},
 ~~\frac{v}{2}-\frac{a}{2c}+\frac{b}{2c},~~\frac{v}{2}+\frac{a}{2c}-\frac{b}{2c}~~~~~~~~~~~~;\\  \frac{v}{2},~1+\frac{v}{2}-\frac{a}{2c}-\frac{b}{2c},~1+\frac{v}{2}+\frac{a}{2c}+\frac{b}{2c},~1+\frac{v}{2}-\frac{a}{2c}+\frac{b}{2c},~1+\frac{v}{2}+\frac{a}{2c}-\frac{b}{2c},~~;\end{array} 1\right),
 \end{multline}
 \end{theorem}
where $\Re(v)<3,~\Re(c)>0,~\Re(vc\pm a\pm b)>0;~\frac{v}{2},~1+\frac{v}{2}\pm\frac{a}{2c}\pm\frac{b}{2c}\in\mathbb{C}\backslash \mathbb{Z}_{0}^{-}$.
\begin{theorem}
The third generalized hyperbolic integral holds true:
 \begin{multline}\label{GRI47}
 \int_{0}^{\infty}\frac{\sinh(ax)\cosh(bx)}{\cosh^{v}(cx)}dx=\frac{2^{v}(v^{2}ac^{2}-a^{3}+ab^{2})}{(vc-a-b)(vc+a+b)(vc-a+b)(vc+a-b)}\times\\\times
 {}_{7}F_{6}\left(\begin{array}{lll}v,1+\sigma_{1},~1+\sigma_{2},~\sigma_{3},~\sigma_{4}
 ,\sigma_{5},~ \sigma_{6}~~~~~~~~~;\\  \sigma_{1},\sigma_{2},1+\sigma_{3} ,~1+\sigma_{4},~1+\sigma_{5},~1+ \sigma_{6};\end{array} -1\right),
 \end{multline}
 \end{theorem}
 where
 \begin{equation*}
 \sigma_{1}= \frac{v}{2}-\frac{\sqrt{a^{2}-b^{2}}}{2c},~\sigma_{2}=\frac{v}{2}+\frac{\sqrt{a^{2}-b^{2}}}{2c}~,
 \sigma_{3}=\frac{v}{2}-\frac{a}{2c}-\frac{b}{2c},
 \end{equation*}
 \begin{equation*}
 ~\sigma_{4}=\frac{v}{2}-\frac{a}{2c}+\frac{b}{2c},~\sigma_{5}=\frac{v}{2}+\frac{a}{2c}+\frac{b}{2c},~
 \sigma_{6}=\frac{v}{2}+\frac{a}{2c}-\frac{b}{2c},
 \end{equation*}
 and $\Re(v)<3,~\Re(c)>0,~\Re(vc\pm a\pm b)>0;~\frac{v}{2}\pm\frac{\sqrt{a^{2}-b^{2}}}{2c},~1+\frac{v}{2}\pm\frac{a}{2c}\pm\frac{b}{2c}\in\mathbb{C}\backslash \mathbb{Z}_{0}^{-}$.
 \begin{theorem}
The fourth generalized hyperbolic integral holds true:
 \begin{multline}\label{GRI48}
 \int_{0}^{\infty}\frac{\sinh(ax)\cosh(bx)}{\sinh^{v}(cx)}dx=\frac{2^{v}(v^{2}ac^{2}-a^{3}+ab^{2})}{(vc-a-b)(vc+a+b)(vc-a+b)(vc+a-b)}\times\\\times
 {}_{7}F_{6}\left(\begin{array}{lll}v,1+\sigma_{1},~1+\sigma_{2},~\sigma_{3},~\sigma_{4}
 ,\sigma_{5},~ \sigma_{6}~~~~~~~~~;\\  \sigma_{1},\sigma_{2},1+\sigma_{3} ,~1+\sigma_{4},~1+\sigma_{5},~1+ \sigma_{6};\end{array} 1\right),
 \end{multline}
 \end{theorem}
 where
 \begin{equation*}
 \sigma_{1}= \frac{v}{2}-\frac{\sqrt{a^{2}-b^{2}}}{2c},~\sigma_{2}=\frac{v}{2}+\frac{\sqrt{a^{2}-b^{2}}}{2c}~,
 \sigma_{3}=\frac{v}{2}-\frac{a}{2c}-\frac{b}{2c},
 \end{equation*}
 \begin{equation*}
 ~\sigma_{4}=\frac{v}{2}-\frac{a}{2c}+\frac{b}{2c},~\sigma_{5}=\frac{v}{2}+\frac{a}{2c}+\frac{b}{2c},~
 \sigma_{6}=\frac{v}{2}+\frac{a}{2c}-\frac{b}{2c},
 \end{equation*}
and $\Re(v)<2,~\Re(c)>0,~\Re(vc\pm a\pm b)>0;~\frac{v}{2}\pm\frac{\sqrt{a^{2}-b^{2}}}{2c},~1+\frac{v}{2}\pm\frac{a}{2c}\pm\frac{b}{2c}\in\mathbb{C}\backslash \mathbb{Z}_{0}^{-}$.
\begin{theorem}
The fifth generalized hyperbolic integral holds true:
 \begin{multline}\label{GRI49}
 \int_{0}^{\infty}\frac{\cosh(ax)\cosh(bx)}{\cosh^{v}(cx)}dx=\frac{2^{v}(v^{3}c^{3}-a^{2}vc-b^{2}vc)}{(vc-a-b)(vc+a+b)(vc-a+b)(vc+a-b)}\times\\\times
 {}_{8}F_{7}\left(\begin{array}{lll}v,~1+\frac{v}{2},~1+\lambda_{1},~1+\lambda_{2},~\sigma_{3},~\sigma_{4}
 ,\sigma_{5},~ \sigma_{6}~;\\ \frac{v}{2},~ \lambda_{1},~\lambda_{2},1+\sigma_{3} ,~1+\sigma_{4},~1+\sigma_{5},~1+\sigma_{6}~;\end{array} -1\right),
 \end{multline}
 \end{theorem}
 where
 \begin{equation*}
 \lambda_{1}= \frac{v}{2}-\frac{\sqrt{a^{2}+b^{2}}}{2c},~\lambda_{2}=\frac{v}{2}+\frac{\sqrt{a^{2}+b^{2}}}{2c}~,
 \sigma_{3}=\frac{v}{2}-\frac{a}{2c}-\frac{b}{2c},
 \end{equation*}
 \begin{equation*}
 ~\sigma_{4}=\frac{v}{2}-\frac{a}{2c}+\frac{b}{2c},~\sigma_{5}=\frac{v}{2}+\frac{a}{2c}+\frac{b}{2c},~
 \sigma_{6}=\frac{v}{2}+\frac{a}{2c}-\frac{b}{2c},
 \end{equation*}
 and $\Re(v)<2,~\Re(c)>0,~\Re(vc\pm a\pm b)>0;\frac{v}{2},~\frac{v}{2}\pm\frac{\sqrt{a^{2}+b^{2}}}{2c},~1+\frac{v}{2}\pm\frac{a}{2c}\pm\frac{b}{2c}\in\mathbb{C}\backslash \mathbb{Z}_{0}^{-}$.
 \begin{theorem}
 The sixth generalized hyperbolic integral holds true:
 \begin{multline}\label{GRI50}
 \int_{0}^{\infty}\frac{\cosh(ax)\cosh(bx)}{\sinh^{v}(cx)}dx=\frac{2^{v}(v^{3}c^{3}-a^{2}vc-b^{2}vc)}{(vc-a-b)(vc+a+b)(vc-a+b)(vc+a-b)}\times\\\times
 {}_{8}F_{7}\left(\begin{array}{lll}v,~1+\frac{v}{2},~1+\lambda_{1},~1+\lambda_{2},~\sigma_{3},~\sigma_{4}
 ,\sigma_{5},~ \sigma_{6}~;\\ \frac{v}{2},~ \lambda_{1},~\lambda_{2},1+\sigma_{3} ,~1+\sigma_{4},~1+\sigma_{5},~1+\sigma_{6}~;\end{array} 1\right),
 \end{multline}
 \end{theorem}
 where
 \begin{equation*}
 \lambda_{1}= \frac{v}{2}-\frac{\sqrt{a^{2}+b^{2}}}{2c},~\lambda_{2}=\frac{v}{2}+\frac{\sqrt{a^{2}+b^{2}}}{2c}~,
 \sigma_{3}=\frac{v}{2}-\frac{a}{2c}-\frac{b}{2c},
 \end{equation*}
 \begin{equation*}
 ~\sigma_{4}=\frac{v}{2}-\frac{a}{2c}+\frac{b}{2c},~\sigma_{5}=\frac{v}{2}+\frac{a}{2c}+\frac{b}{2c},~
 \sigma_{6}=\frac{v}{2}+\frac{a}{2c}-\frac{b}{2c},
 \end{equation*}
 and $\Re(v)<1,~\Re(c)>0,~\Re(vc\pm a\pm b)>0;~\frac{v}{2},~\frac{v}{2}\pm\frac{\sqrt{a^{2}+b^{2}}}{2c},~1+\frac{v}{2}\pm\frac{a}{2c}\pm\frac{b}{2c}\in\mathbb{C}\backslash \mathbb{Z}_{0}^{-}$.\\
\textbf{Hypergeometric proof of integral (\ref{GRI49})}:   Suppose left hand side of eq.(\ref{GRI49}) is denoted by $\Upsilon(a,b,c,v)$ and using the product formula of hyperbolic function in the left hand side of eq.(\ref{GRI49}), we get
   \begin{equation}\label{GRI78}
  \Upsilon(a,b,c,v)= \frac{1}{2}\int_{0}^{\infty}\frac{\cosh\{(a+b)x\}}{\cosh^{v}(c x)}dx+\frac{1}{2}\int_{0}^{\infty}\frac{\cosh\{(a-b)x\}}{\cosh^{v}(c x)}dx=\textbf{L}_{1}+\textbf{L}_{2},
 \end{equation}
 where $\textbf{L}_{1}$ and  $\textbf{L}_{2}$ are given by
 \begin{equation}\label{GRI79}
 \textbf{L}_{1}= \frac{1}{2}\int_{0}^{\infty}\frac{\cosh\{(a+b)x\}}{\cosh^{v}(c x)}dx ~~and~~\textbf{L}_{2}=\frac{1}{2}\int_{0}^{\infty}\frac{\cosh\{(a-b)x\}}{\cosh^{v}(c x)}dx.
 \end{equation}
 Using exponential definition of hyperbolic functions in the integral $\textbf{L}_{1}$, which yields
 \begin{equation*}\label{GRI80}
\textbf{L}_{1}= 2^{v-2}\int_{0}^{\infty}e^{-vcx}\bigg[e^{(a+b)x}+e^{-(a+b)x}\bigg]\left(1+e^{-2cx}\right)^{-v}dx.
 \end{equation*}
 \begin{equation}\label{GRI81}
= 2^{v-2}\int_{0}^{\infty}e^{-vcx}\bigg[e^{(a+b)x}+e^{-(a+b)x}\bigg] {}_{1}F_{0}\left(\begin{array}{lll}v~;\\  \overline{~~~};\end{array} -e^{-2cx} \right)dx,
 \end{equation}
  when $\Re(c)>0$ , then $|-e^{-2cx}|<1$ for all real $x>0$. It is the convergence condition of above binomial function  ${}_{1}F_{0}(\cdot)$ in eq.(\ref{GRI81}), then it yields
  \begin{equation}\label{GRI82}
 \textbf{L}_{1} =2^{v-2}\sum_{r=0}^{\infty}\frac{(v)_{r}}{r!}(-1)^{r}\bigg[\int_{0}^{\infty}e^{-\{(vc-a-b)+2cr\}x}dx+\int_{0}^{\infty}e^{-\{(vc+a+b)+2cr\}x}dx\bigg],
 \end{equation}
 where $\Re(vc-a-b)>0,~\Re(vc+a+b)>0,~\Re(c)>0$, it is the convergence conditions of Laplace transform of unity in the integral (\ref{GRI82}). Then
applying Laplace transformation formula (\ref{GRI22}) in the eq.(\ref{GRI82}), we obtain
 \begin{equation}\label{GRI83}
 \textbf{L}_{1}=2^{v-2}\sum_{r=0}^{\infty}\frac{(v)_{r}}{r!}(-1)^{r}\bigg[\frac{1}{\{(vc-a-b)+2cr\}}+\frac{1}{\{(vc+a+b)+2cr\}}\bigg],
 \end{equation}
~~~~ where $\Re(vc-a-b)>0,~\Re(vc+a+b)>0,~\Re(c)>0$.\\
 Similarly, proof of $\textbf{L}_{2}$ is given by
 \begin{equation}\label{GRI84}
  \textbf{L}_{2}=2^{v-2}\sum_{r=0}^{\infty}\frac{(v)_{r}}{r!}(-1)^{r}\bigg[\frac{1}{\{(vc-a+b)+2cr\}}+\frac{1}{\{(vc+a-b)+2cr\}}\bigg],
 \end{equation}
 ~~~~ where $\Re(vc-a+b)>0,~\Re(vc+a-b)>0,~\Re(c)>0$.\\
  Making use of the eqns (\ref{GRI83}) and (\ref{GRI84}) in the above eq. (\ref{GRI78}), we obtain
  \begin{multline*}\label{GRI85}
 \Upsilon(a,b,c,v)=2^{v-2}\sum_{r=0}^{\infty}\frac{(v)_{r}}{r!}(-1)^{r}\times\\\times \bigg[\frac{1}{\{(vc-a-b)+2cr\}}+\frac{1}{\{(vc+a+b)+2cr\}}
 +\frac{1}{\{(vc-a+b)+2cr\}}+\frac{1}{\{(vc+a-b)+2cr\}}\bigg],
 \end{multline*}
  \begin{eqnarray}\label{GRI86}
 =2^{v-1}\sum_{r=0}^{\infty}\frac{(v)_{r}}{r!}(-1)^{r}\times\nonumber~~~~~~~~~~~~~~~~~~~~~~~~~~~~~~~~~~~~~~~~~~~~~~~~~~~~~~~~~~~~~~~~~~~~~~~~~~~~~~~~~~~~~~~~~~~\\\times
 \bigg[\frac{(vc-a+2cr)}{\{(vc-a-b)+2cr\}\{(vc-a+b)+2cr\}}+\frac{(vc+a+2cr)}{\{(vc+a+b)+2cr\}\{(vc+a-b)+2cr\}}\bigg]\nonumber\\
 \end{eqnarray}
 ~~~~ where $\Re(vc\pm a\pm b)>0,~\Re(c)>0$.
 After simplifications we obtain
  \begin{eqnarray*}
 \Upsilon(a,b,c,v)=2^{v-1}\sum_{r=0}^{\infty}\frac{(v)_{r}}{r!}(-1)^{r}\times\nonumber~~~~~~~~~~~~~~~~~~~~~~~~~~~~~~~~~~~~~~~~~~~~~~~~~~~~~~~~~~~~~~~~~~~~~~~~~~~~~~~~~~~~~~~~~~~\\\times
 \bigg[\frac{16c^{3}r^{3}+24vc^{3}r^{2}+(12v^{2}c^{3}-4a^{2}c-4b^{2}c)r + (2v^{3}c^{3}-2va^{2}c-2b^{2}vc)}{\{(vc-a-b)+2cr\}\{(vc+a+b)+2cr\}\{(vc-a+b)+2cr\}\{(vc+a-b)+2cr\}}\bigg],\nonumber\\
 \end{eqnarray*}
 \begin{eqnarray}\label{GRI87}
 =2^{v-1}\sum_{r=0}^{\infty}\frac{(v)_{r}}{r!}(-1)^{r}\times\nonumber~~~~~~~~~~~~~~~~~~~~~~~~~~~~~~~~~~~~~~~~~~~~~~~~~~~~~~~~~~~~~~~~~~~~~~~~~~~~~~~~~~~~~~~~~~~\\\times
 \bigg[\frac{16c^{3}\left\{r-\left(\frac{-vc+\sqrt{a^{2}+b^{2}}}{2c}\right)\right\}\left\{r-\left(\frac{-vc-\sqrt{a^{2}+b^{2}}}{2c}\right)\right\}\left\{r+\frac{v}{2}\right\}}
 {\{(vc-a-b)+2cr\}\{(vc+a+b)+2cr\}\{(vc-a+b)+2cr\}\{(vc+a-b)+2cr\}}\bigg],\nonumber\\
 \end{eqnarray}
 ~~~~ where $\Re(vc\pm a\pm b)>0,~\Re(c)>0$.\\
  Employ algebraic properties of Pochhammer symbol in the eq.(\ref{GRI87}), after simplifications, we obtain
 \begin{eqnarray*}\label{GRI89}
 \Upsilon(a,b,c,v)=\frac{2^{v}(v^{3}c^{3}-a^{2}vc-b^{2}vc)}{(vc-a-b)(vc+a+b)(vc-a+b)(vc+a-b)}\sum_{r=0}^{\infty}\bigg[\frac{(v)_{r}\left(1+\frac{v}{2}\right)_{r}}{\left(\frac{v}{2}\right)_{r}r!}\times~~~~~~~~~~~~~\\ \times\frac{\left(1+\frac{v}{2}-\frac{\sqrt{a^{2}+b^{2}}}{2c}\right)_{r}\left(1+\frac{v}{2}+\frac{\sqrt{a^{2}+b^{2}}}{2c}\right)_{r}\left(\frac{vc-a-b}{2c}\right)_{r}\left(\frac{vc+a+b}{2c}
 \right)_{r}\left(\frac{vc-a+b}{2c}\right)_{r}\left(\frac{vc+a-b}{2c}\right)_{r}(-1)^{r}}
 {\left(\frac{v}{2}-\frac{\sqrt{a^{2}+b^{2}}}{2c}\right)_{r}\left(\frac{v}{2}+\frac{\sqrt{a^{2}+b^{2}}}{2c}\right)_{r}\left(\frac{vc-a-b+2c}{2c}\right)_{r}\left(\frac{vc+a+b+2c}{2c}\right)_{r}\left(\frac{vc-a+b+2c}{2c}\right)_{r}\left(\frac{vc+a-b+2c}{2c}\right)_{r}}\bigg],
 \end{eqnarray*}
 \begin{eqnarray}\label{GRI90}
 =\frac{2^{v}(v^{3}c^{3}-a^{2}vc-b^{2}vc)}{(vc-a-b)(vc+a+b)(vc-a+b)(vc+a-b)}\times~~~~~~~~~~~~~~~~~~~~~~~~\nonumber\\ \times
 {}_{8}F_{7}\left(\begin{array}{lll}v,~1+\frac{v}{2},~1+\lambda_{1},~1+\lambda_{2},~\sigma_{3},~\sigma_{4}
 ,\sigma_{5},~ \sigma_{6}~;\\ \frac{v}{2},~ \lambda_{1},~\lambda_{2},1+\sigma_{3} ,~1+\sigma_{4},~1+\sigma_{5},~1+\sigma_{6}~;\end{array} -1\right),
 \end{eqnarray}
 where
 \begin{equation*}
 \lambda_{1}= \frac{v}{2}-\frac{\sqrt{a^{2}+b^{2}}}{2c},~\lambda_{2}=\frac{v}{2}+\frac{\sqrt{a^{2}+b^{2}}}{2c}~,
 \sigma_{3}=\frac{v}{2}-\frac{a}{2c}-\frac{b}{2c},
 \end{equation*}
 \begin{equation*}
 ~\sigma_{4}=\frac{v}{2}-\frac{a}{2c}+\frac{b}{2c},~\sigma_{5}=\frac{v}{2}+\frac{a}{2c}+\frac{b}{2c},~
 \sigma_{6}=\frac{v}{2}+\frac{a}{2c}-\frac{b}{2c},
 \end{equation*}
and $\Re(c)>0,~\Re(vc\pm a\pm b)>0;\frac{v}{2},~\frac{v}{2}\pm\frac{\sqrt{a^{2}+b^{2}}}{2c},~1+\frac{v}{2}\pm\frac{a}{2c}\pm\frac{b}{2c}\in\mathbb{C}\backslash \mathbb{Z}_{0}^{-}$.
  Similarly, proofs of the integrals (\ref{GRI45}),(\ref{GRI46}),(\ref{GRI47}),(\ref{GRI48}) and (\ref{GRI50}) are much akin to that of the integral (\ref{GRI49}), which we have already discussed in a detailed manner.
\section{Special class of some hyperbolic integrals in terms of ${}_{3}F_{2}(\pm1),{}_{4}F_{3}(\pm1)$ and ${}_{5}F_{4}(\pm1)$}
\texttt{Each of the following hyperbolic definite integrals holds true:}

 $\bullet$ When $v=2$ and $c=b$ in the eq.(\ref{GRI45}), we get
 \begin{equation}\label{FA3}
 \int_{0}^{\infty}\frac{\sinh(ax)\sinh(bx)}{\cosh^{2}(bx)}dx=\frac{16ab^{2}}{(b^{2}-a^{2})(9b^{2}-a^{2})}
 ~~{}_{4}F_{3}\left(\begin{array}{lll}2,2,~\frac{1}{2}+\frac{a}{2b},~\frac{1}{2}-\frac{a}{2b};\\  1,~\frac{5}{2}+\frac{a}{2b},~\frac{5}{2}-\frac{a}{2b}~~~;\end{array} -1\right),
 \end{equation}
 \begin{equation}\label{FA00}
 ~~~=\frac{a\pi}{2b^{2}}\sec\left(\frac{\pi a}{2b}\right),~~~~~~~~~~~~~~~~~~~~~~~~~~~~~~~~~~~~~~~~
 \end{equation}
where $\Re(b)>0,~\Re(b\pm a)>0,~\Re(3b\pm a)>0,~\frac{5}{2}\pm\frac{a}{2b}\in\mathbb{C}\backslash \mathbb{Z}_{0}^{-}$;~ $\frac{a}{b}\in\mathbb{C}\backslash\{\pm1,\pm3,\pm5,...\}$.
Using hypergeometric form of $\sec (z) $  function (\ref{GRI32}) [when $z=\frac{\pi a}{2b}$ ]  in the eq. (\ref{FA3}), we obtain right hand side of (\ref{FA00}). Also, right hand side of eq.(\ref{FA00}) can be obtained by using summation theorem (\ref{GRI12}) and recurrence relation for gamma function in the hypergeometric series (\ref{FA3}).\\
 $\bullet$ When $v=2$ and $c=b$ in the eq.(\ref{GRI46}), we get
 \begin{equation}\label{FA33}
 \int_{0}^{\infty}\frac{\sinh(ax)}{\sinh(bx)}dx=\frac{16ab^{2}}{(b^{2}-a^{2})(9b^{2}-a^{2})}
 ~~{}_{4}F_{3}\left(\begin{array}{lll}2,2,~\frac{1}{2}+\frac{a}{2b},~\frac{1}{2}-\frac{a}{2b};\\  1,~\frac{5}{2}+\frac{a}{2b},~\frac{5}{2}-\frac{a}{2b}~~~;\end{array} 1\right),
 \end{equation}
 \begin{equation}\label{FA000}
 ~~~=\frac{\pi}{2b}\tan\left(\frac{\pi a}{2b}\right),~~~~~~~~~~~~~~~~~~~~~~~~~~~~~~~~~~~~~~~~
 \end{equation}
 where $\Re(b)>0,~\Re(b\pm a)>0,~\Re(3b\pm a)>0,~\frac{5}{2}\pm\frac{a}{2b}\in\mathbb{C}\backslash \mathbb{Z}_{0}^{-}$;~ $\frac{a}{b}\in\mathbb{C}\backslash\{\pm1,\pm3,\pm5,...\}$.
The right hand side of eq.(\ref{FA000}) can be obtained by using summation theorem (\ref{GRI102}), recurrence relations for gamma function in the hypergeometric series (\ref{FA33}).\\
$\bullet$ When $v=2$, $b=a$, $c=1$ in the eq.(\ref{GRI46}), we get
 \begin{eqnarray}\label{FA0011}
 \frac{1}{2}\int_{-\infty}^{\infty}\frac{\sinh^{2}(a x)}{\sinh^{2}( x)}dx=\int_{0}^{\infty}\frac{\sinh^{2}(a x)}{\sinh^{2}( x)}dx
 =\frac{a^{2}}{(1-a^{2})}~~{}_{3}F_{2}\left(\begin{array}{lll}1,~1-a,~1+a;\\ 2-a,~2+a~~~~;\end{array} 1\right),
 \end{eqnarray}
 \begin{equation}
 =\frac{a}{2}[\Psi(1+a)-\Psi(1-a)],
 \end{equation}
 \begin{equation}\label{FA0012}
 =\frac{1}{2}[1-a\pi\cot(a\pi)],~~~~~~~~~~~
 \end{equation}
 where $a\neq0,\pm1,\pm2,\pm3,...$\\
The right hand side of eq.(\ref{FA0012}) can be obtained by using the properties of Digamma function (\ref{GRI24})-(\ref{GRI262}) in the hypergeometric series (\ref{FA0011}).\\
$\bullet$  In the eq.(\ref{GRI47})  we interchange $a$ and $b$; $v=2$, then put $c=b$, we get
 \begin{multline}\label{FA001}
 \int_{0}^{\infty}\frac{\cosh(ax)\sinh(bx)}{\cosh^{2}(bx)}dx\\=\frac{(12b^{3}+4a^{2}b)}{(b^{2}-a^{2})(9b^{2}-a^{2})}
 ~~{}_{5}F_{4}\left(\begin{array}{lll}2,~2-\frac{\sqrt{b^{2}-a^{2}}}{2b},~2+\frac{\sqrt{b^{2}-a^{2}}}{2b},~\frac{1}{2}-\frac{a}{2b},~\frac{1}{2}+\frac{a}{2b};\\  1-\frac{\sqrt{b^{2}-a^{2}}}{2b},~1+\frac{\sqrt{b^{2}-a^{2}}}{2b},~\frac{5}{2}+\frac{a}{2b},~\frac{5}{2}-\frac{a}{2b}~~~~;\end{array} -1\right),
 \end{multline}
where $\Re(b)>0,~\Re(b\pm a)>0,~\Re(3b\pm a)>0;~1\pm\frac{\sqrt{b^{2}-a^{2}}}{2b},~\frac{5}{2}\pm\frac{a}{2b}\in\mathbb{C}\backslash \mathbb{Z}_{0}^{-}$ .\\
$\bullet$ When  $b=0$ in the eq.(\ref{GRI47}), we get
 \begin{equation}\label{cc1}
 \int_{0}^{\infty}\frac{\sinh(ax)}{\cosh^{v}(cx)}dx=\frac{2^{v}a}{[(vc)^{2}-a^{2}]}
 ~~{}_{3}F_{2}\left(\begin{array}{lll}v,~\frac{v}{2}-\frac{a}{2c},~\frac{v}{2}+\frac{a}{2c};\\  1+\frac{v}{2}-\frac{a}{2c},~1+\frac{v}{2}+\frac{a}{2c};\end{array} -1\right),
 \end{equation}
 \begin{equation}\label{cc7}
= \frac{2^{v-1}}{(vc-a)}~{}_{2}F_{1}\left(\begin{array}{lll}v,~\frac{v}{2}-\frac{a}{2c}~~;\\  1+\frac{v}{2}-\frac{a}{2c};\end{array} -1\right)
 -\frac{2^{v-1}}{(vc+a)}~{}_{2}F_{1}\left(\begin{array}{lll}v,~\frac{v}{2}+\frac{a}{2c}~~;\\  1+\frac{v}{2}+\frac{a}{2c};\end{array} -1\right),
 \end{equation}
 where $\Re(v)<3,~\Re(c)>0,~\Re(vc\pm a)>0,~~1+\frac{v}{2}\pm\frac{a}{2c}\in\mathbb{C}\backslash \mathbb{Z}_{0}^{-}$.\\
 $\bullet$ When  $b=0$ in the eq.(\ref{GRI48}), we get
 \begin{equation}\label{cc2}
 \int_{0}^{\infty}\frac{\sinh(ax)}{\sinh^{v}(cx)}dx=\frac{2^{v}a}{[(vc)^{2}-a^{2}]}
 ~~{}_{3}F_{2}\left(\begin{array}{lll}v,~\frac{v}{2}-\frac{a}{2c},~\frac{v}{2}+\frac{a}{2c};\\  1+\frac{v}{2}-\frac{a}{2c},~1+\frac{v}{2}+\frac{a}{2c};\end{array} 1\right),
 \end{equation}
  \begin{equation}\label{cc2}
 =\frac{2^{v-2}}{(c)\Gamma(v)}\Gamma\left(\frac{v}{2}+\frac{a}{2c}\right)\Gamma\left(\frac{v}{2}-\frac{a}{2c}\right)
 \frac{\sin\left(\frac{a\pi}{2c}\right)}{\sin\left(\frac{v\pi}{2}\right)},
 \end{equation}
  where $\Re(v)<2,~\Re(c)>0,~\Re(vc\pm a)>0,~~1+\frac{v}{2}\pm\frac{a}{2c}\in\mathbb{C}\backslash \mathbb{Z}_{0}^{-}$ and $v\neq0,\pm2,\pm4,\pm6,...$\\
$\bullet$ When $v=1$ and $b=0$, then $c=b$ in the eq.(\ref{GRI48}), we get
 \begin{eqnarray}\label{FA1}
 \int_{0}^{\infty}\frac{\sinh(ax)}{\sinh(b x)}dx= \frac{2a}{(b^{2}-a^{2})}~~{}_{3}F_{2}\left(\begin{array}{lll}1,~\frac{1}{2}-\frac{a}{2b},~~\frac{1}{2}+\frac{a}{2b};\\  \frac{3}{2}-\frac{a}{2b},~~\frac{3}{2}+\frac{a}{2b}~~~~;\end{array} 1\right),
 \end{eqnarray}
 \begin{equation}\label{FA49}
=\frac{\pi}{2b}\tan\left(\frac{\pi a}{2b}\right),~~~~~~~~~~~~~~~~~~~~~~~~~~~~
 \end{equation}
where $\Re(b)>0,~\Re(b\pm a)>0,~\frac{3}{2}\pm\frac{a}{2b}\in\mathbb{C}\backslash \mathbb{Z}_{0}^{-}$;~$\frac{a}{b}\in\mathbb{C}\backslash\{\pm1,\pm3,\pm5,...\}$.\\
Using hypergeometric form of $\tan (z) $  function (\ref{GRI31}) [ when $z=\frac{\pi a}{2b}$ ]  in the eq. (\ref{FA1}), we obtain right hand side of eq.(\ref{FA49}). Also, right hand side of eq.(\ref{FA49}) can be obtained  by using the Dixon's theorem ${}_{3}F_{2}(1)$ (\ref{GRI11}) in the eq. (\ref{FA1}).\\
 $\bullet$ When  $b=0$ in the eq.(\ref{GRI49}), we get
 \begin{equation}
 \int_{0}^{\infty}\frac{\cosh(ax)}{\cosh^{v}(cx)}dx=\frac{2^{v}(vc)}{[(vc)^{2}-a^{2}]}
 ~~{}_{4}F_{3}\left(\begin{array}{lll}v,~1+\frac{v}{2},\frac{v}{2}-\frac{a}{2c},~\frac{v}{2}+\frac{a}{2c};\\  \frac{v}{2}, 1+\frac{v}{2}-\frac{a}{2c},~1+\frac{v}{2}+\frac{a}{2c}~~~;\end{array} -1\right),
 \end{equation}
  \begin{equation}\label{cc3}
 =\frac{2^{v-2}}{(c)\Gamma(v)}\Gamma\left(\frac{v}{2}+\frac{a}{2c}\right)\Gamma\left(\frac{v}{2}-\frac{a}{2c}\right),
 \end{equation}
  where $\Re(v)<2,~\Re(c)>0,~\Re(vc\pm a)>0;~~\frac{v}{2},1+\frac{v}{2}\pm\frac{a}{2c}\in\mathbb{C}\backslash \mathbb{Z}_{0}^{-}$.\\
$\bullet$ When $v=1$ and $b=0$, then $c=b$ in the eq.(\ref{GRI49}), we get
 \begin{eqnarray}\label{FA2}
 \int_{0}^{\infty}\frac{\cosh(ax)}{\cosh(bx)}dx=\frac{2b}{(b^{2}-a^{2})}~~{}_{4}F_{3}\left(\begin{array}{lll}1,~\frac{3}{2},~\frac{1}{2}-\frac{a}{2b},~~\frac{1}{2}+\frac{a}{2b};\\  \frac{1}{2},~\frac{3}{2}-\frac{a}{2b},~~\frac{3}{2}+\frac{a}{2b}~~~~;\end{array} -1\right),
 \end{eqnarray}
 \begin{equation}\label{FA50}
 =\frac{\pi}{2b}\sec\left(\frac{\pi a}{2b}\right),~~~~~~~~~~~~~~~~~~~~~~~~~~~~~~~~~~~
 \end{equation}
where $\Re(b)>0,~\Re(b\pm a)>0,~\frac{3}{2}\pm\frac{a}{2b}\in\mathbb{C}\backslash \mathbb{Z}_{0}^{-}$;~ $\frac{a}{b}\in\mathbb{C}\backslash\{\pm1,\pm3,\pm5,...\}$.\\
Using hypergeometric form of $\sec (z) $  function (\ref{GRI32}) [when $z=\frac{\pi a}{2b}$ ]  in the eq. (\ref{FA2}), we obtain right hand side of eq. (\ref{FA50}).  Also, right hand side of eq.(\ref{FA50}) can be obtained  by using the classical summation theorem ${}_{4}F_{3}(-1)$ (\ref{GRI12}) in the eq. (\ref{FA2}).\\
 $\bullet$ In the eq.(\ref{FA50}) replacing $a\rightarrow 2ia $ and $b=\pi$, we get  a known result of Ramanujan \cite[p.11,eq.(1.5.1(27))]{E1}
 \begin{equation}\label{F00}
  \int_{0}^{\infty}\frac{\cos(2ax)}{\cosh(\pi x)}dx
 =\frac{1}{2}sech\left(a\right),~~~~~~~~~|\Im(a)|<\frac{\pi}{2}.
 \end{equation}
 $\bullet$ When  $b=0$ in the eq.(\ref{GRI50}), we get
 \begin{equation}
 \int_{0}^{\infty}\frac{\cosh(ax)}{\sinh^{v}(cx)}dx=\frac{2^{v}(vc)}{[(vc)^{2}-a^{2}]}
 ~~{}_{4}F_{3}\left(\begin{array}{lll}v,~1+\frac{v}{2},\frac{v}{2}-\frac{a}{2c},~\frac{v}{2}+\frac{a}{2c};\\  \frac{v}{2}, 1+\frac{v}{2}-\frac{a}{2c},~1+\frac{v}{2}+\frac{a}{2c}~~~;\end{array} 1\right),
 \end{equation}
  \begin{equation}\label{cc4}
 =\frac{2^{v-2}}{(c)\Gamma(v)}\Gamma\left(\frac{v}{2}+\frac{a}{2c}\right)\Gamma\left(\frac{v}{2}-\frac{a}{2c}\right)
 \frac{\cos\left(\frac{a\pi}{2c}\right)}{\cos\left(\frac{v\pi}{2}\right)},
 \end{equation}
 where $\Re(v)<1,~\Re(c)>0,~\Re(vc\pm a)>0;~~\frac{v}{2},1+\frac{v}{2}\pm\frac{a}{2c}\in\mathbb{C}\backslash \mathbb{Z}_{0}^{-}$ and $v\neq\pm1,\pm3,\pm5,...$\\
$\bullet$ When $v=1$ and $b=0$, then $c=b$ in the eq.(\ref{GRI47}), we get
 \begin{equation}\label{FA54}
\int_{0}^{\infty}\frac{\sinh(ax)}{\cosh(b x)}dx =\frac{2a}{(b^{2}-a^{2})}{}_{3}F_{2}\left(\begin{array}{lll}1,~\frac{1}{2}-\frac{a}{2b},~~\frac{1}{2}+\frac{a}{2b};\\  \frac{3}{2}-\frac{a}{2b},~~\frac{3}{2}+\frac{a}{2b}~~~~;\end{array} -1\right),
 \end{equation}
 \begin{equation}\label{FA503}
~~~=\frac{\pi}{2b}\sec\left(\frac{\pi a}{2b}\right)-\frac{1}{b}\beta\left(\frac{a+b}{2b}\right),
\end{equation}
\begin{equation}\label{FA53}
~~~~~~~~~~~~~~~~~~~~~~~~~~~~~~~~~=\frac{\pi}{2b}\sec\left(\frac{\pi a}{2b}\right)-\frac{1}{2b}\left[\Psi\left(\frac{a+3b}{4b}\right)-\Psi\left(\frac{a+b}{4b}\right)\right],
 \end{equation}
 \begin{equation}\label{F00}
=\frac{(b^{2}-a^{2})}{8ab}\left[\Psi\left(\frac{3b-a}{4b}\right)-\Psi\left(\frac{b-a}{4b}\right)-\Psi\left(\frac{3b+a}{4b}\right)+\Psi\left(\frac{b+a}{4b}\right)\right],
 \end{equation}
where $\Re(b)>0,~\Re(b\pm a)>0,~\frac{3}{2}\pm\frac{a}{2b}\in\mathbb{C}\backslash \mathbb{Z}_{0}^{-}$;~ $\frac{a}{b}\in\mathbb{C}\backslash\{\pm1,\pm3,\pm5,...\}$ .\\
  \textbf{Independent proofs of} (\ref{FA54})-(\ref{F00}): Taking left hand side of eq.(\ref{FA54}) and suppose it is denoted by $\Phi(a,b)$ upon using the well known result of hyperbolic function , we get
  \begin{equation}\label{AF04}
 \Phi(a,b)
 =\int_{0}^{\infty}\left(\frac{e^{ax}}{e^{bx}+e^{-bx}}\right)dx-\int_{0}^{\infty}\left(\frac{e^{-ax}}{e^{bx}+e^{-bx}}\right)dx=\textbf{Y}_{1}-\textbf{Y}_{2}, \end{equation}
where $\textbf{Y}_{1}$ and $\textbf{Y}_{2}$ are given by
 \begin{equation}\label{AF05}
 \textbf{Y}_{1}= \int_{0}^{\infty}\left(\frac{e^{ax}}{e^{bx}+e^{-bx}}\right)dx,~~and~~ \textbf{Y}_{2}= \int_{0}^{\infty}\left(\frac{e^{-ax}}{e^{bx}+e^{-bx}}\right)dx.
 \end{equation}
 From the above integral $\textbf{Y}_{1}$ can also written by
\begin{equation*}
\textbf{Y}_{1}= \int_{0}^{\infty}e^{-(b-a)x}\bigg(1+e^{-2bx}\bigg)^{-1}dx,
\end{equation*}
 \begin{equation}\label{AF06}
= \int_{0}^{\infty}\bigg[e^{-(b-a)x}{}_{1}F_{0}\left(\begin{array}{lll}1~;\\  \overline{~~~};\end{array} -e^{-2bx} \right)\bigg]dx,
 \end{equation}
 when $\Re(b)>0$ , then $|-e^{-2bx}|<1$ for all $x>0$. It is the convergence conditions of above binomial function  ${}_{1}F_{0}(\cdot)$ in eq.(\ref{AF06}), then it yields
 \begin{equation*}
 \textbf{Y}_{1}=\int_{0}^{\infty}\bigg[e^{-(b-a)x}\sum_{r=0}^{\infty}(-1)^{r}~e^{-2brx}\bigg]dx,
 \end{equation*}
  \begin{equation}\label{AF07}
  =\sum_{r=0}^{\infty}(-1)^{r}\bigg[\int_{0}^{\infty}e^{-\{(b-a)+2br\}x}dx\bigg],
 \end{equation}
 where $\Re(b-a)>0,~\Re(b)>0$, it is the convergence condition of Laplace transform of unity in the integral (\ref{AF07}). Then
applying formula (\ref{GRI22}) in the eq.(\ref{AF07}), we obtain
 \begin{equation}\label{AF08}
 \textbf{Y}_{1}=\sum_{r=0}^{\infty}(-1)^{r}\bigg[\frac{1}{(b-a)+2br}\bigg], ~~~~~~~~~~~~~\Re(b-a)>0,~\Re(b)>0.
 \end{equation}
 Similarly, proof of $\textbf{Y}_{2}$ is given by
 \begin{equation}\label{AF09}
 \textbf{Y}_{2}=\sum_{r=0}^{\infty}(-1)^{r}\bigg[\frac{1}{(b+a)+2br}\bigg], ~~~~~~~~~~~~~\Re(b+a)>0,~\Re(b)>0.
 \end{equation}
 Making use of the eqns (\ref{AF08}) and (\ref{AF09}) in the above eq. (\ref{AF04}), we obtain
 \begin{equation}\label{AF10}
\Phi(a,b)=\sum_{r=0}^{\infty}\bigg[\frac{(-1)^{r}}{(b-a)+2br}-\frac{(-1)^{r}}{(b+a)+2br}\bigg],
 \end{equation}
 ~~~~~where $\Re(b\pm a)>0,~\Re(b)>0.$\\
 Employ algebraic properties of Pochhammer symbol in the eq.(\ref{AF10}), after simplifications, we obtain
 \begin{equation*}\label{AF11}
 \Phi(a,b)=\frac{2a}{b^{2}-a^{2}}\sum_{r=0}^{\infty}\bigg[\frac{\left(\frac{b-a}{2b}\right)_{r}\left(\frac{b+a}{2b}\right)_{r}(-1)^{r}}
 {\left(\frac{3b-a}{2b}\right)_{r}\left(\frac{3b+a}{2b}\right)_{r}}\bigg],
 \end{equation*}
 \begin{equation}\label{AF12}
 ~~~~~~~~~~=\frac{2a}{(b^{2}-a^{2})}~~{}_{3}F_{2}\left(\begin{array}{lll}1,~\frac{1}{2}-\frac{a}{2b},~~\frac{1}{2}+\frac{a}{2b};\\  \frac{3}{2}-\frac{a}{2b},~~\frac{3}{2}+\frac{a}{2b};\end{array} -1\right),
 \end{equation}
 where $\Re(b\pm a)>0,~\Re(b)>0,~\frac{a}{b}\in\mathbb{C}\backslash\{\pm1,\pm3,\pm5,...\},~\frac{3}{2}\pm\frac{a}{2b}\in\mathbb{C}\backslash \mathbb{Z}_{0}^{-}$. Now, proof of the eq.(\ref{FA503}) is obtained by using the eq. (\ref{AF10}) with addition and substraction of its second term, is given below
 \begin{equation}\label{AF101}
\Phi(a,b)=\sum_{r=0}^{\infty}\bigg[\frac{(-1)^{r}}{(b-a)+2br}+\frac{(-1)^{r}}{(b+a)+2br}-\frac{2(-1)^{r}}{(b+a)+2br}\bigg],
 \end{equation}
 \begin{equation}\label{AF102}
=\sum_{r=0}^{\infty}(-1)^{r}\bigg[\frac{(1+2r)}{\{(b-a)+2br\}\{(b+a)+2br\}}\bigg]-\frac{1}{b}\sum_{r=0}^{\infty}\frac{(-1)^{r}}{\big\{\left(\frac{b+a}{2b}\right)+r\big\}}.
 \end{equation}
 Employ algebraic properties of Pochhammer symbol and Lower case beta function (\ref{GRI25}) in the eq.(\ref{AF102}), after simplifications, we obtain
 \begin{equation*}\label{AF103}
 \Phi(a,b)=\frac{2b}{b^{2}-a^{2}}\sum_{r=0}^{\infty}\bigg[\frac{\left(\frac{3}{2}\right)_{r}\left(\frac{b-a}{2b}\right)_{r}\left(\frac{b+a}{2b}\right)_{r}(-1)^{r}}
 {\left(\frac{1}{2}\right)_{r}\left(\frac{3b-a}{2b}\right)_{r}\left(\frac{3b+a}{2b}\right)_{r}}\bigg]-\frac{1}{b}\beta\left(\frac{a+b}{2b}\right),
 \end{equation*}
 \begin{equation}\label{AF104}
 ~~~~~~~~~~=\frac{2b}{(b^{2}-a^{2})}~~{}_{4}F_{3}\left(\begin{array}{lll}1,~\frac{3}{2},~\frac{1}{2}-\frac{a}{2b},~~\frac{1}{2}+\frac{a}{2b};\\  \frac{1}{2},~\frac{3}{2}-\frac{a}{2b},~~\frac{3}{2}+\frac{a}{2b};\end{array} -1\right)-\frac{1}{b}\beta\left(\frac{a+b}{2b}\right),
 \end{equation}
where $\Re(b\pm a)>0,~\Re(b)>0,~\frac{a}{b}\in\mathbb{C}\backslash\{\pm1,\pm3,\pm5,...\},~\frac{3}{2}\pm\frac{a}{2b}\in\mathbb{C}\backslash \mathbb{Z}_{0}^{-}$. We obtain, upon using hypergeometric function of sec(z) function (\ref{GRI32}) [when $z=\frac{\pi a}{2b}$ ] in the eq. (\ref{AF104}).
\begin{equation}\label{AF1014}
\Phi(a,b)=\frac{\pi}{2b}\sec\left(\frac{\pi a}{2b}\right)-\frac{1}{b}\beta\left(\frac{a+b}{2b}\right).
\end{equation}
Using the classical summation theorem ${}_{4}F_{3}(-1)$ (\ref{GRI12}) in the right hand side of eq.(\ref{AF104}), after simplification we get eq.(\ref{AF1014}) or (\ref{FA503}). In the eq.(\ref{FA54}) apply the properties of Digamma functions (\ref{GRI26}) and (\ref{GRI000}), we get the result (\ref{F00}).
  \section{Applications of product formulas in special class of hyperbolic integrals}
   Product formulas of hyperbolic functions:
 \begin{equation}\label{GRI34}
 \sinh(A)\cosh(B)=\frac{1}{2}[\sinh(A+B)+\sinh(A-B)],
 \end{equation}
 \begin{equation}\label{GRI35}
 \sinh(A)\sinh(B)=\frac{1}{2}[\cosh(A+B)-\cosh(A-B)],
 \end{equation}
 \begin{equation}\label{GRI36}
 \cosh(A)\cosh(B)=\frac{1}{2}[\cosh(A+B)+\cosh(A-B)].
 \end{equation}
\texttt{Each of the following hyperbolic definite integrals holds true:}\\
$\bullet$ In the eq.(\ref{GRI45}) using product formula and applying the result (\ref{cc3}), we get
\begin{multline}\label{sss1}
   \int_{0}^{\infty}\frac{\sinh(ax)\sinh(bx)}{\cosh^{v}(c x)}dx\\
   =\frac{2^{v-3}}{(c)\Gamma(v)}\bigg[\Gamma\left(\frac{vc+a+b}{2c}\right)\Gamma\left(\frac{vc-a-b}{2c}\right)
   -\Gamma\left(\frac{vc+a-b}{2c}\right)\Gamma\left(\frac{vc-a+b}{2c}\right)\bigg],
\end{multline}
where $\Re(v)<4,~\Re(c)>0,~\Re(vc\pm a\pm b)>0;~\frac{v}{2},~1+\frac{v}{2}\pm\frac{a}{2c}\pm\frac{b}{2c}\in\mathbb{C}\backslash \mathbb{Z}_{0}^{-}$.\\
 $\bullet$ In the eq.(\ref{GRI45}) put $v=1$ and interchange $b$ and $c$, we get
\begin{multline}\label{FAF811}
   \int_{0}^{\infty}\frac{\sinh(ax)\sinh(cx)}{\cosh(b x)}dx
   =\frac{2^{2}acb}{(b-a-c)(b+a+c)(b-a+c)(b+a-c)}\times\\\times
 {}_{6}F_{5}\left(\begin{array}{lll}1,\frac{3}{2},\frac{1}{2}-\frac{a}{2b}-\frac{c}{2b},
 \frac{1}{2}-\frac{a}{2b}+\frac{c}{2b},\frac{1}{2}+\frac{a}{2b}+\frac{c}{2b},\frac{1}{2}+\frac{a}{2b}-\frac{c}{2b}~;\\  \frac{1}{2},\frac{3}{2}-\frac{a}{2b}-\frac{c}{2b}
 ,\frac{3}{2}-\frac{a}{2b}+\frac{c}{2b},\frac{3}{2}+\frac{a}{2b}+\frac{c}{2b},\frac{3}{2}+\frac{a}{2b}-\frac{c}{2b};\end{array} -1\right),
 \end{multline}
 where $\Re(b)>0,~\Re(b\pm a\pm c)>0,~\frac{3}{2}\pm\frac{a}{2b}\pm\frac{c}{2b}\in\mathbb{C}\backslash \mathbb{Z}_{0}^{-}$.\\
 $\bullet$  Using the product formula of hyperbolic function, then applying eq.(\ref{FA50}), we get
 \begin{eqnarray}\label{FAF8}
    \int_{0}^{\infty}\frac{\sinh(ax)\sinh(cx)}{\cosh(bx)}dx
    =\frac{\pi}{4b}\left[\sec\left(\frac{\pi(a+c)}{2b}\right)-\sec\left(\frac{\pi(a-c)}{2b}\right)\right],
 \end{eqnarray}
 \begin{equation}\label{FAF81}
 =\left(\frac{\pi}{b}\right)\frac{\sin\left(\frac{a\pi}{2b}\right)\sin\left(\frac{c\pi}{2b}\right)}
 {\left\{\cos\left(\frac{c\pi}{b}\right)+\cos\left(\frac{a\pi}{b}\right)\right\}},
 \end{equation}
  where $\Re(b)>0,~\Re(b\pm a\pm c)>0,~\frac{3}{2}\pm\frac{a}{2b}\pm\frac{c}{2b}\in\mathbb{C}\backslash \mathbb{Z}_{0}^{-}$.\\
 $\bullet$ In the eq.(\ref{GRI46}) using product formula and applying the result (\ref{cc4}), we get
\begin{multline}\label{ss2}
   \int_{0}^{\infty}\frac{\sinh(ax)\sinh(bx)}{\sinh^{v}(c x)}dx
   =\frac{2^{v-3}}{(c)\Gamma(v)}\bigg[\Gamma\left(\frac{vc+a+b}{2c}\right)\Gamma\left(\frac{vc-a-b}{2c}\right)
   \frac{\cos\left(\frac{(a+b)\pi}{2c}\right)}{\cos(\frac{v\pi}{2})}\\
   -\Gamma\left(\frac{vc+a-b}{2c}\right)\Gamma\left(\frac{vc-a+b}{2c}\right)\frac{\cos\left(\frac{(a-b)\pi}{2c}\right)}{\cos(\frac{v\pi}{2})}\bigg],
\end{multline}
where $\Re(v)<3,~\Re(c)>0,~\Re(vc\pm a\pm b)>0;~\frac{v}{2},~1+\frac{v}{2}\pm\frac{a}{2c}\pm\frac{b}{2c}\in\mathbb{C}\backslash \mathbb{Z}_{0}^{-}$.\\
 $\bullet$ In the eq.(\ref{GRI47}) using product formula and then applying the result (\ref{cc7}), we get
 \begin{multline}\label{cc9}
   \int_{0}^{\infty}\frac{\sinh(ax)\cosh(bx)}{\cosh^{v}(c x)}dx\\
=2^{v-2} \bigg[\frac{1}{(vc-a-b)}~{}_{2}F_{1}\left(\begin{array}{lll}v,\frac{vc-a-b}{2c}~~~;\\ 1+\frac{vc-a-b}{2c};\end{array} -1\right)
 -\frac{1}{(vc+a+b)}~{}_{2}F_{1}\left(\begin{array}{lll}v,\frac{vc+a+b}{2c}~~~;\\ 1+\frac{vc+a+b}{2c};\end{array} -1\right)\\
 +\frac{1}{(vc-a+b)}~{}_{2}F_{1}\left(\begin{array}{lll}v,\frac{vc-a+b}{2c}~~~;\\ 1+\frac{vc-a+b}{2c};\end{array} -1\right)
 -\frac{1}{(vc+a-b)}~{}_{2}F_{1}\left(\begin{array}{lll}v,\frac{vc+a-b}{2c}~~~;\\ 1+\frac{vc+a-b}{2c};\end{array} -1\right)\bigg],
 \end{multline}
where  $\Re(v)<3,~\Re(c)>0,~\Re(vc\pm a\pm b)>0;~\frac{v}{2}\pm\frac{\sqrt{a^{2}-b^{2}}}{2c},~1+\frac{v}{2}\pm\frac{a}{2c}\pm\frac{b}{2c}\in\mathbb{C}\backslash \mathbb{Z}_{0}^{-}$.\\
 $\bullet$ In the eq.(\ref{GRI47}) put $v=1$ and interchange $b$ and $c$, we get
\begin{multline}\label{FAF111}
   \int_{0}^{\infty}\frac{\sinh(ax)\cosh(cx)}{\cosh(b x)}dx
   =\frac{2(ab^{2}-a^{3}+ac^{2})}{(b-a-c)(b+a+c)(b-a+c)(b+a-c)}\times\\\times
 {}_{7}F_{6}\left(\begin{array}{lll}1,\frac{3}{2}-\frac{\sqrt{a^{2}-c^{2}}}{2b},\frac{3}{2}+\frac{\sqrt{a^{2}-c^{2}}}{2b},\frac{1}{2}-\frac{a}{2b}-\frac{c}{2b},
 \frac{1}{2}-\frac{a}{2b}+\frac{c}{2b},\frac{1}{2}+\frac{a}{2b}+\frac{c}{2b},\frac{1}{2}+\frac{a}{2b}-\frac{c}{2b}~;\\  \frac{1}{2}-\frac{\sqrt{a^{2}-c^{2}}}{2b},\frac{1}{2}+\frac{\sqrt{a^{2}-c^{2}}}{2b},\frac{3}{2}-\frac{a}{2b}-\frac{c}{2b}
 ,\frac{3}{2}-\frac{a}{2b}+\frac{c}{2b},\frac{3}{2}+\frac{a}{2b}+\frac{c}{2b},\frac{3}{2}+\frac{a}{2b}-\frac{c}{2b};\end{array} -1\right),
 \end{multline}
 where $\Re(b)>0,~\Re(b\pm a\pm c)>0;~\frac{1}{2}\pm\frac{\sqrt{a^{2}-c^{2}}}{2b},~\frac{3}{2}\pm\frac{a}{2b}\pm\frac{c}{2b}\in\mathbb{C}\backslash \mathbb{Z}_{0}^{-}$.\\
 $\bullet$   Using the product formula of hyperbolic function , then applying eq.(\ref{FA503}), we get
 \begin{eqnarray}\label{FAF10}
    \int_{0}^{\infty}\frac{\sinh(ax)\cosh(cx)}{\cosh(bx)}dx
    =\frac{\pi}{4b}\left[\sec\left(\frac{\pi(a+c)}{2b}\right)+\sec\left(\frac{\pi(a-c)}{2b}\right)\right]\nonumber\\
    -\frac{1}{2b}\left[\beta\left(\frac{(a+b+c)}{2b}\right)+\beta\left(\frac{(a+b-c)}{2b}\right)\right],
 \end{eqnarray}
 \begin{equation}\label{FAF101}
     =\left(\frac{\pi}{b}\right)\frac{\cos\left(\frac{a\pi}{2b}\right)\cos\left(\frac{c\pi}{2b}\right)}
     {\left\{\cos\left(\frac{c\pi}{b}\right)+\cos\left(\frac{a\pi}{b}\right)\right\}}
     -\frac{1}{2b}\left[\beta\left(\frac{(a+b+c)}{2b}\right)+\beta\left(\frac{(a+b-c)}{2b}\right)\right],
 \end{equation}
  where $\Re(b)>0,~\Re(b\pm a\pm c)>0;~\frac{1}{2}\pm\frac{\sqrt{a^{2}-c^{2}}}{2b},~\frac{3}{2}\pm\frac{a}{2b}\pm\frac{c}{2b}\in\mathbb{C}\backslash \mathbb{Z}_{0}^{-}$.\\
$\bullet$ In the eq.(\ref{GRI48}) using product formula and applying the result (\ref{cc2}), we get
\begin{multline}\label{ss3}
   \int_{0}^{\infty}\frac{\sinh(ax)\cosh(bx)}{\sinh^{v}(c x)}dx
   =\frac{2^{v-3}}{(c)\Gamma(v)}\bigg[\Gamma\left(\frac{vc+a+b}{2c}\right)\Gamma\left(\frac{vc-a-b}{2c}\right)
   \frac{\sin\left(\frac{(a+b)\pi}{2c}\right)}{\sin(\frac{v\pi}{2})}\\
   +\Gamma\left(\frac{vc+a-b}{2c}\right)\Gamma\left(\frac{vc-a+b}{2c}\right)\frac{\sin\left(\frac{(a-b)\pi}{2c}\right)}{\sin(\frac{v\pi}{2})}\bigg],
\end{multline}
where $\Re(v)<2,~\Re(c)>0,~\Re(vc\pm a\pm b)>0;~\frac{v}{2}\pm\frac{\sqrt{a^{2}-b^{2}}}{2c},~1+\frac{v}{2}\pm\frac{a}{2c}\pm\frac{b}{2c}\in\mathbb{C}\backslash \mathbb{Z}_{0}^{-}$.\\
$\bullet$ In the eq.(\ref{GRI48}) put $v=1$ and interchange $b$ and $c$, we get
\begin{multline}\label{FAF710}
   \int_{0}^{\infty}\frac{\sinh(ax)\cosh(cx)}{\sinh(b x)}dx
   =\frac{2(ab^{2}-a^{3}+ac^{2})}{(b-a-c)(b+a+c)(b-a+c)(b+a-c)}\times\\\times
 {}_{7}F_{6}\left(\begin{array}{lll}1,\frac{3}{2}-\frac{\sqrt{a^{2}-c^{2}}}{2b},\frac{3}{2}+\frac{\sqrt{a^{2}-c^{2}}}{2b},\frac{1}{2}-\frac{a}{2b}-\frac{c}{2b},
 \frac{1}{2}-\frac{a}{2b}+\frac{c}{2b},\frac{1}{2}+\frac{a}{2b}+\frac{c}{2b},\frac{1}{2}+\frac{a}{2b}-\frac{c}{2b}~;\\  \frac{1}{2}-\frac{\sqrt{a^{2}-c^{2}}}{2b},\frac{1}{2}+\frac{\sqrt{a^{2}-c^{2}}}{2b},\frac{3}{2}-\frac{a}{2b}-\frac{c}{2b}
 ,\frac{3}{2}-\frac{a}{2b}+\frac{c}{2b},\frac{3}{2}+\frac{a}{2b}+\frac{c}{2b},\frac{3}{2}+\frac{a}{2b}-\frac{c}{2b};\end{array} 1\right),
 \end{multline}
 where $\Re(b)>0,~\Re(b\pm a\pm c)>0;~\frac{1}{2}\pm\frac{\sqrt{a^{2}-c^{2}}}{2b},~\frac{3}{2}\pm\frac{a}{2b}\pm\frac{c}{2b}\in\mathbb{C}\backslash \mathbb{Z}_{0}^{-}$.\\
  $\bullet$  Using the product formula of hyperbolic function , then applying eq.(\ref{FA49}), we get
\begin{eqnarray}\label{FAF7}
    \int_{0}^{\infty}\frac{\sinh(ax)\cosh(cx)}{\sinh(b x)}dx
    =\frac{\pi}{4b}\left[\tan\left(\frac{\pi(a+c)}{2b}\right)+\tan\left(\frac{\pi(a-c)}{2b}\right)\right],
 \end{eqnarray}
 \begin{equation}\label{FAF71}
    =\left(\frac{\pi}{2b}\right)\frac{\sin\left(\frac{a\pi}{b}\right)}
    {\left\{\cos\left(\frac{c\pi}{b}\right)+\cos\left(\frac{a\pi}{b}\right)\right\}},~
 \end{equation}
  where $\Re(b)>0,~\Re(b\pm a\pm c)>0;~\frac{1}{2}\pm\frac{\sqrt{a^{2}-c^{2}}}{2b},~\frac{3}{2}\pm\frac{a}{2b}\pm\frac{c}{2b}\in\mathbb{C}\backslash \mathbb{Z}_{0}^{-}$.\\
 ~
$\bullet$ In the eq.(\ref{GRI49}) using product formula and applying the result (\ref{cc3}), we get
\begin{multline}\label{ss4}
   \int_{0}^{\infty}\frac{\cosh(ax)\cosh(bx)}{\cosh^{v}(c x)}dx\\
   =\frac{2^{v-3}}{(c)\Gamma(v)}\bigg[\Gamma\left(\frac{vc+a+b}{2c}\right)\Gamma\left(\frac{vc-a-b}{2c}\right)
   +\Gamma\left(\frac{vc+a-b}{2c}\right)\Gamma\left(\frac{vc-a+b}{2c}\right)\bigg],
\end{multline}
where $\Re(v)<2,~\Re(c)>0,~\Re(vc\pm a\pm b)>0;\frac{v}{2},~\frac{v}{2}\pm\frac{\sqrt{a^{2}+b^{2}}}{2c},~1+\frac{v}{2}\pm\frac{a}{2c}\pm\frac{b}{2c}\in\mathbb{C}\backslash \mathbb{Z}_{0}^{-}$.\\
$\bullet$ In the eq.(\ref{GRI49}) put $v=1$ and interchange $b$ and $c$, we get
\begin{multline}\label{FAF109}
   \int_{0}^{\infty}\frac{\cosh(ax)\cosh(cx)}{\cosh(b x)}dx
   =\frac{2(b^{3}-a^{2}b-c^{2}b)}{(b-a-c)(b+a+c)(b-a+c)(b+a-c)}\times\\\times
 {}_{8}F_{7}\left(\begin{array}{lll}1,\frac{3}{2},\frac{3}{2}-\frac{\sqrt{a^{2}+c^{2}}}{2b},\frac{3}{2}+\frac{\sqrt{a^{2}+c^{2}}}{2b},\frac{1}{2}-\frac{a}{2b}-\frac{c}{2b},
 \frac{1}{2}-\frac{a}{2b}+\frac{c}{2b},\frac{1}{2}+\frac{a}{2b}+\frac{c}{2b},\frac{1}{2}+\frac{a}{2b}-\frac{c}{2b}~;\\  \frac{1}{2}, \frac{1}{2}-\frac{\sqrt{a^{2}+c^{2}}}{2b},\frac{1}{2}+\frac{\sqrt{a^{2}+c^{2}}}{2b},\frac{3}{2}-\frac{a}{2b}-\frac{c}{2b}
 ,\frac{3}{2}-\frac{a}{2b}+\frac{c}{2b},\frac{3}{2}+\frac{a}{2b}+\frac{c}{2b},\frac{3}{2}+\frac{a}{2b}-\frac{c}{2b};\end{array} -1\right),
 \end{multline}
 where $\Re(b)>0,~\Re(b\pm a\pm c)>0;~\frac{1}{2}\pm\frac{\sqrt{a^{2}+c^{2}}}{2b},~\frac{3}{2}\pm\frac{a}{2b}\pm\frac{c}{2b}\in\mathbb{C}\backslash \mathbb{Z}_{0}^{-}$.\\
 $\bullet$  Using the product formula of hyperbolic function, then applying eq.(\ref{FA50}), we get
 \begin{eqnarray}\label{FAF9}
    \int_{0}^{\infty}\frac{\cosh(ax)\cosh(cx)}{\cosh(bx)}dx
    =\frac{\pi}{4b}\left[\sec\left(\frac{\pi(a+c)}{2b}\right)+\sec\left(\frac{\pi(a-c)}{2b}\right)\right],
    \end{eqnarray}
    \begin{equation}\label{FAF09}
  =\left(\frac{\pi}{b}\right)\frac{\cos\left(\frac{a\pi}{2b}\right)\cos\left(\frac{c\pi}{2b}\right)}
  {\left\{\cos\left(\frac{c\pi}{b}\right)+\cos\left(\frac{a\pi}{b}\right)\right\}},
 \end{equation}
  where $\Re(b)>0,~\Re(b\pm a\pm c)>0;~\frac{1}{2}\pm\frac{\sqrt{a^{2}+c^{2}}}{2b},~\frac{3}{2}\pm\frac{a}{2b}\pm\frac{c}{2b}\in\mathbb{C}\backslash \mathbb{Z}_{0}^{-}$.\\
 $\bullet$ In the eq.(\ref{GRI50}) using product formula and applying the result (\ref{cc4}), we get
\begin{multline}\label{ss5}
   \int_{0}^{\infty}\frac{\cosh(ax)\cosh(bx)}{\sinh^{v}(c x)}dx
   =\frac{2^{v-3}}{(c)\Gamma(v)}\bigg[\Gamma\left(\frac{vc+a+b}{2c}\right)\Gamma\left(\frac{vc-a-b}{2c}\right)
   \frac{\cos\left(\frac{(a+b)\pi}{2c}\right)}{\cos(\frac{v\pi}{2})}\\
   +\Gamma\left(\frac{vc+a-b}{2c}\right)\Gamma\left(\frac{vc-a+b}{2c}\right)\frac{\cos\left(\frac{(a-b)\pi}{2c}\right)}{\cos(\frac{v\pi}{2})}\bigg],
\end{multline}
where $\Re(v)<1,~\Re(c)>0,~\Re(vc\pm a\pm b)>0;~\frac{v}{2},~\frac{v}{2}\pm\frac{\sqrt{a^{2}+b^{2}}}{2c},~1+\frac{v}{2}\pm\frac{a}{2c}\pm\frac{b}{2c}\in\mathbb{C}\backslash \mathbb{Z}_{0}^{-}$.
\section*{Conclusion}
Here, we have described some definite integrals containing the quotients of hyperbolic functions. Thus certain integrals of hyperbolic functions, which may be different from those of presented here, can also be evaluated in a similar way. Therefore, the results presented in this paper can be expressed in terms of hypergeometric functions, trigonometric and hyperbolic functions, Digamma functions, Beta function of one variable, and Gamma function.

\end{document}